\documentclass[12pt]{article}

\usepackage[T1]{fontenc}
\usepackage[utf8]{inputenc}
\usepackage[margin=1in]{geometry}
\usepackage{amsmath,amssymb,amsthm,mathtools,mathrsfs}
\usepackage{times}
\usepackage{needspace}
\usepackage[
  colorlinks=true,
  linkcolor=blue,
  citecolor=blue,
  urlcolor=blue
]{hyperref}

\hypersetup{
  pdftitle={Inverse energy transfer for 2D Navier--Stokes near a self-similar vortex collapse},
  pdfauthor={Trinh T. Nguyen}
}

\allowdisplaybreaks
\numberwithin{equation}{section}

\newtheorem*{theorem}{Theorem}
\newtheorem{corollary}{Corollary}
\newtheorem{proposition}{Proposition}[section]
\newtheorem{lemma}[proposition]{Lemma}
\theoremstyle{remark}

\newcommand{\dd}{\,\mathrm d}

\newcommand{\R}{\mathbb R}
\newcommand{\wh}{\widehat}
\newcommand{\ol}{\overline}
\newcommand{\cH}{\mathcal H}
\newcommand{\cM}{\mathcal M}
\newcommand{\cR}{\mathcal R}
\newcommand{\sE}{\mathscr E}
\newcommand{\Gtot}{\Gamma_{\mathrm{tot}}}
\newcommand{\Gabs}{\Gamma_{\mathrm{abs}}}
\newcommand{\ind}{\mathbf1}
\newcommand{\Real}{\operatorname{Re}}

\title{Inverse energy transfer for 2D Navier--Stokes\\
near a self-similar vortex collapse}

\author{Trinh T. Nguyen\thanks{
School of Mathematical Sciences,
Shanghai Jiao Tong University, Shanghai, China.\newline
Email: \texttt{q0004349@sjtu.edu.cn}.
}}

\date{}

\begin{document}
\maketitle

\begin{abstract}
\noindent
We prove inverse energy transfer for the two-dimensional Navier--Stokes
equations near a self-similar three-vortex collapse, confirming the 1979
prediction of Novikov and Sedov. Nonlinear interactions transfer kinetic
energy from high to low frequencies over a common time
interval that shrinks to the collapse time as the viscosity tends to zero.
We prove that the net nonlinear transfer over this interval converges to the same positive constant across an expanding range of frequency cutoffs, with quantitative error bounds uniform in the cutoff.
\end{abstract}

\setcounter{tocdepth}{1}
\tableofcontents

\section{Introduction}\label{sec:introduction}

We study inverse energy transfer near a self-similar three-vortex
collapse in the two-dimensional Navier--Stokes equations as the viscosity
tends to zero. Nonlinear interactions redistribute kinetic energy among
spatial scales in turbulent flows~\cite{CET94,CNT99,Ey05,ES06,BV19,DE19,CDP26}.
In two dimensions, the inverse energy cascade carries kinetic energy
toward lower wavenumbers and larger spatial scales
\cite{Kr67,Co97,Ey06,Drivas19,BCZPSW20,CDP25,Drivas26}. In domains with boundaries, spectral and nonlinear
instabilities of boundary layers illustrate the formation of small
spatial scales~\cite{GGN16,GN18,GN19,GN24}.

It is therefore natural to ask whether this transfer can be quantified
for specific flow configurations. Vortex dynamics provides such configurations. Its mathematical study
goes back to Helmholtz and Kelvin~\cite{He58,Th68}.
Localized concentrations of vorticity can persist over long times and
organize the surrounding flow through their motion and interaction.
In two dimensions, point-vortex systems give a finite-dimensional
description of this motion~\cite{Ne01}.
The stability, dynamics, and inviscid limits of vortex flows
have been studied in
\cite{MP94,GW05,Ga11,Ga12,GSV19,NN19,GSV24,DG26,DKK26,NguyenVortexWave,NguyenCollapse,NguyenAmplification}.
Point-vortex dynamics can also be derived from the Boltzmann
equation~\cite{KN25}. Vortex interactions and the evolution of energy
and enstrophy have been studied in inviscid, viscous, and regularized
models~\cite{McW84,MZM88,MLL05,Sa12,GS16,GS18,Go25,DT26,GSV26}.

For inviscid point vortices, Novikov computed the point-vortex energy spectrum and related its
evolution to the mean wavenumber~\cite{No75}.
Novikov and Sedov applied this calculation to a self-similar
three-vortex collapse and interpreted the decrease of the interaction
moment as inverse energy transfer~\cite{NS79}.
We recall their argument and derive the finite-cutoff formulas in
Section~\ref{sec:novikov-sedov}.

In this paper, we consider the two-dimensional Navier--Stokes equations on
$\R^2$ in vorticity form,
\begin{equation}\label{eq:NS}
\partial_t\omega^\nu+u^\nu\cdot\nabla\omega^\nu
=\nu\Delta\omega^\nu,
\qquad u^\nu=K*\omega^\nu,
\qquad K(x)=\frac{x^\perp}{2\pi|x|^2},
\end{equation}
where $\nu>0$, $x^\perp=(-x_2,x_1)$, and the initial vorticity is
\begin{equation}\label{eq:initial-data}
\omega^\nu(0)=\sum_{k=1}^3\Gamma_k\delta_{z_k^0}.
\end{equation}
The strengths $\Gamma_1,\Gamma_2,\Gamma_3$ are fixed nonzero real
numbers, and the positions $z_1^0,z_2^0,z_3^0$ are distinct. All are
independent of $\nu$. For each $\nu>0$, the initial datum determines a
unique global solution $\omega^\nu$, smooth at positive times~\cite{GG05,Ga11}.

We choose a configuration whose point-vortex evolution collapses
self-similarly at time $1$~\cite{Ar10,NguyenCollapse}. We assume that,
for some $\beta\in\R$, the initial configuration satisfies
\begin{equation}\label{eq:configuration}
\frac12z_k^0-\beta(z_k^0)^\perp
+\sum_{j\ne k}\Gamma_jK(z_k^0-z_j^0)=0,
\qquad k=1,2,3.
\end{equation}
The associated point-vortex trajectories are
\begin{equation}\label{eq:centers}
z_k(t)=\sqrt{1-t}\,
\mathfrak R_{-\beta\log(1-t)}z_k^0,
\qquad 0\le t<1,
\end{equation}
where $\mathfrak R_\theta$ denotes counterclockwise rotation
through the angle $\theta$. They solve
$\dot z_k=\sum_{j\ne k}\Gamma_jK(z_k-z_j)$ and satisfy
\begin{equation}\label{eq:distances}
r_{ij}(t)=|z_i(t)-z_j(t)|=\sqrt{1-t}\,r_{ij}(0),
\qquad r_{ij}(0)=|z_i^0-z_j^0|.
\end{equation}
We write $r_{\min}(0)=\min_{i<j}r_{ij}(0)$.

The point-vortex Hamiltonian is
\begin{equation}\label{eq:transfer-constant}
\cH=-\frac1{2\pi}\sum_{i<j}\Gamma_i\Gamma_j\log r_{ij}(0)<0.
\end{equation}
We prove that the net energy transferred from higher to lower frequencies
converges to $-\cH>0$ as $\nu\to0$. We prescribe the time remaining until collapse and the frequency
cutoffs as powers of the viscosity,
with the exponents fixed before taking the limit. An explicit error
bound holds uniformly over this frequency range and hence across a
number of dyadic cutoffs proportional to $|\log\nu|$. The limiting
transfer is the same across all these cutoffs during a common interval
whose length is asymptotic to a prescribed power of $\nu$.

\paragraph{Notation.}
We set $\Gtot=\sum_k\Gamma_k$ and
$\Gabs=\sum_k|\Gamma_k|$.
We denote the indicator of a set $D$ by $\ind_D$ and use the Fourier convention
$\wh f(\xi)=\int_{\R^2}e^{-ix\cdot\xi}f(x)\dd x$.
In Fourier integrals, we write $\wh f(t)$ for
$\wh f(t,\xi)$ when the integration variable is clear;
$\Real$ and an overline denote the real part and
complex conjugation, respectively.
For nonnegative quantities $a,b$, we write $a\lesssim b$ if
$a\le Cb$ for a positive constant $C$.
Unless otherwise specified, implicit constants may depend on the fixed
initial configuration but are independent of the viscosity, time, and
frequency cutoffs. In the proof of the main theorem, they may also
depend on the fixed exponents $\alpha_\pm,\gamma_\pm$.
Subscripts indicate additional dependence.
For $b>0$, the notation $a=O(b)$, $a=o(b)$, and $a\sim b$
means, respectively, $|a|\lesssim b$, $a/b\to0$, and $a/b\to\text{const}$
in the limit under consideration.

\section{Main results}\label{sec:main-results}

\begin{theorem}
We fix $0<\gamma_-<\alpha_-<\alpha_+<\gamma_+<1$ and set $A_\nu=1-\nu^{\gamma_-}$ and $B_\nu=1-\nu^{\gamma_+}$. There exist constants $C>0$ and $\nu_0\in(0,1)$, depending
only on the fixed initial configuration and the exponents
$\alpha_\pm,\gamma_\pm$, such that, for every $0<\nu\le\nu_0$
and every cutoff $\nu^{-\alpha_-/2}\le N\le\nu^{-\alpha_+/2}$, the following estimates hold:
\begin{equation}\label{eq:plateau-estimates}
\begin{aligned}
&\left|
-\frac1{(2\pi)^2}\int_{A_\nu}^{B_\nu}\Real
\int_{\R^2}\ind_{\{|\xi|>N\}}
\ol{\wh{u^\nu}(t)}\cdot
\wh{(u^\nu\cdot\nabla)u^\nu}(t)\dd\xi\dd t
-\cH\right|\\
&\quad\le C\left(
\nu^{1-\gamma_+}(1+|\log\nu|)+\nu^{\gamma_-}
+\nu^{2(\gamma_+-\alpha_+)}
+\nu^{3(\alpha_--\gamma_-)/4}\right),\\[6pt]
&\left|
\frac1{2(2\pi)^2}\int_{\R^2}\ind_{\{|\xi|>N\}}
\Bigl(|\wh{u^\nu}(B_\nu)|^2
-|\wh{u^\nu}(A_\nu)|^2\Bigr)\dd\xi
-\cH\right|\\
&\quad\le C\left(
\nu^{1-\gamma_+}(1+|\log\nu|)+\nu^{\gamma_-}
+\nu^{2(\gamma_+-\alpha_+)}
+\nu^{3(\alpha_--\gamma_-)/4}\right).
\end{aligned}
\end{equation}
In particular, these bounds hold uniformly over the dyadic cutoffs
\begin{equation}\label{eq:dyadic-cutoffs}
N_{j,\nu}=2^j\nu^{-\alpha_-/2},\qquad 0\le j\le J_\nu,\qquad J_\nu=\left\lfloor
\frac{\alpha_+-\alpha_-}{2\log2}|\log\nu|
\right\rfloor\longrightarrow\infty.
\end{equation}
The net transfer toward lower frequencies is $-\cH+o(1)$
across each of $J_\nu+1$ distinct cutoffs during the same interval
$[A_\nu,B_\nu]$, with an error uniform in $j$.
\end{theorem}

For $0<\nu<1$, the observation times satisfy
$0<A_\nu<B_\nu<1$, and the length of the common interval tends
to zero:
\begin{equation}\label{eq:intro-duration}
B_\nu-A_\nu
\sim\nu^{\gamma_-}\longrightarrow0.
\end{equation}
During the same interval $[A_\nu,B_\nu]$, energy is transferred toward
lower frequencies across an increasing number of cutoffs as $\nu\to0$.
The dyadic cutoffs are
\[
N_{0,\nu},\qquad 2N_{0,\nu},\qquad 4N_{0,\nu},\qquad
\ldots,\qquad 2^{J_\nu}N_{0,\nu},
\]
where $N_{0,\nu}=\nu^{-\alpha_-/2}$.
For each cutoff $N_{j,\nu}$, we measure the net nonlinear energy
transferred from $|\xi|>N_{j,\nu}$ to $|\xi|\le N_{j,\nu}$ over the
common interval $[A_\nu,B_\nu]$. Every cutoff records
$-\cH+o(1)$, with the error uniform in $j$. Thus the same positive
limiting transfer is measured across a number of cutoffs that grows
like $|\log\nu|$.
The following corollary gives the energy gain below each cutoff.

\begin{corollary}\label{cor:energy-gain}
Under the assumptions of the main theorem, we fix
$k_0>0$. As $\nu\to0$, the following limits hold uniformly
for $\nu^{-\alpha_-/2}\le N\le\nu^{-\alpha_+/2}$:
\begin{equation}\label{eq:band-gain}
\begin{split}
&
\frac1{2(2\pi)^2}\int_{\R^2}\ind_{\{k_0<|\xi|\le N\}}
\Bigl(|\wh{u^\nu}(B_\nu)|^2
-|\wh{u^\nu}(A_\nu)|^2\Bigr)\dd\xi
\longrightarrow-\cH,\\[6pt]
&
-\frac1{(2\pi)^2}\int_{A_\nu}^{B_\nu}\Real
\int_{\R^2}\ind_{\{k_0<|\xi|\le N\}}
\ol{\wh{u^\nu}(t)}\cdot
\wh{(u^\nu\cdot\nabla)u^\nu}(t)\dd\xi\dd t
\longrightarrow-\cH.
\end{split}
\end{equation}
\end{corollary}

\section{Ideas of the proof}\label{sec:proof-ideas}

We compute the nonlinear transfer from the change in kinetic energy
above a Fourier cutoff. For each viscosity, the cutoff $N$ is held
fixed throughout the observation interval $[A_\nu,B_\nu]$.
Integrating the energy identity \eqref{eq:high-balance} gives
\[
\begin{aligned}
&\frac1{(2\pi)^2}\int_{A_\nu}^{B_\nu}\Real
\int_{\R^2}\ind_{\{|\xi|>N\}}
\ol{\wh u^\nu(t)}\cdot
\wh{(u^\nu\cdot\nabla)u^\nu}(t)\dd\xi\dd t\\
&\overset{\eqref{eq:high-balance}}{=}
\underbrace{\frac1{2(2\pi)^2}\int_{\R^2}\ind_{\{|\xi|>N\}}
\bigl(|\wh u^\nu(A_\nu)|^2-|\wh u^\nu(B_\nu)|^2\bigr)\dd\xi}_{(\mathrm{I})}
-\underbrace{\frac{\nu}{(2\pi)^2}
\int_{A_\nu}^{B_\nu}\int_{\R^2}\ind_{\{|\xi|>N\}}
|\wh\omega^\nu(t)|^2\dd\xi\dd t}_{(\mathrm{II})}.
\end{aligned}
\]
Here $(\mathrm I)$ is the decrease of energy above $N$, and
$(\mathrm{II})$ is the energy lost there by viscosity. Their difference
is the net nonlinear transfer toward lower frequencies. Our two goals
are therefore to prove $(\mathrm I)\to-\cH>0$ and
$(\mathrm{II})\to0$, uniformly for
$\nu^{-\alpha_-/2}\le N\le\nu^{-\alpha_+/2}$.

The viscous loss is controlled directly by the vorticity smoothing
estimate. Corollary~\ref{cor:vanishing-viscous-dissipation} gives
\[
0\le(\mathrm{II})
\le\nu\int_{A_\nu}^{B_\nu}\|\omega^\nu(t)\|_2^2\dd t
\lesssim\Gabs^2\log\frac{B_\nu}{A_\nu}
\lesssim\nu^{\gamma_-}.
\]
This bound is independent of $N$. It tends to zero because the
observation interval shrinks around the positive time $1$.
Thus an order-one decrease of high-frequency energy on this interval
must be accounted for by the nonlinear transfer.

We estimate $(\mathrm I)$ by comparing $u^\nu$ with $u_H^\nu=K*Q^\nu$,
where $Q^\nu$ is the sum of Gaussians centered at the point vortices
$z_j(t)$, defined in \eqref{eq:Gaussian}. The sharp approximation
bound of \cite{NguyenCollapse}, recalled in
Proposition~\ref{prop:approximation}, gives
\[
\sup_{A_\nu\le t\le B_\nu}
r_\nu(t)\lesssim\nu^{1-\gamma_+},
\qquad
r_\nu(t)=\|\omega^\nu(t)-Q^\nu(t)\|_{L^1(\R^2)}.
\]
It controls the Gaussian approximation up to
$B_\nu=1-\nu^{\gamma_+}$, only $\nu^{\gamma_+}$ before collapse.
Applying the logarithmic energy estimate \eqref{eq:high-energy-error}
in Lemma~\ref{lem:high-energy-error} at $A_\nu$ and $B_\nu$, we obtain,
uniformly over the prescribed cutoffs,
\[
\begin{aligned}
(\mathrm I)
&=\frac1{2(2\pi)^2}\int_{\R^2}\ind_{\{|\xi|>N\}}
\bigl(|\wh u_H^\nu(A_\nu)|^2-|\wh u_H^\nu(B_\nu)|^2\bigr)\dd\xi\\
&\qquad+O\!\left(\nu^{1-\gamma_+}(1+|\log\nu|)\right).
\end{aligned}
\]

We next compute the asymptotics of the Gaussian energy explicitly.
The spectral representation in Lemma~\ref{lem:Gaussian-high} and the
interaction estimates in Lemma~\ref{lem:interaction-estimates} give,
uniformly over the prescribed cutoffs,
\[
\begin{aligned}
&\frac1{2(2\pi)^2}\int_{\R^2}\ind_{\{|\xi|>N\}}
\bigl(|\wh u_H^\nu(A_\nu)|^2-|\wh u_H^\nu(B_\nu)|^2\bigr)\dd\xi\\
&\qquad=-\cH+O\!\left(\nu^{\gamma_-}
+\nu^{2(\gamma_+-\alpha_+)}
+\nu^{3(\alpha_--\gamma_-)/4}\right).
\end{aligned}
\]
Combining this expansion with the comparison estimate above and the
bound on $(\mathrm{II})$, we obtain
\[
\begin{aligned}
(\mathrm I)-(\mathrm{II})
&=-\cH+O\!\left(\nu^{1-\gamma_+}(1+|\log\nu|)
+\nu^{\gamma_-}+\nu^{2(\gamma_+-\alpha_+)}
+\nu^{3(\alpha_--\gamma_-)/4}\right)\\
&\longrightarrow-\cH>0,
\end{aligned}
\]
where the strict positivity follows from
Proposition~\ref{thm:spectral-positivity}. 
\section{Novikov--Sedov theory}
\label{sec:novikov-sedov}

We recall the point-vortex spectrum and the inverse-transfer prediction of Novikov and Sedov~\cite{No75,NS79}. We derive finite-cutoff formulas for the mean wavenumber. The Bessel representation also enters the Navier--Stokes energy estimates.

\subsection{The point-vortex spectrum}

Angular Fourier averages can be expressed in terms of vortex distances through the Bessel function of the first kind of order zero. We define this function by
\begin{equation}\label{eq:bessel-definition}
J_0(r):=\frac1{2\pi}\int_0^{2\pi}\cos(r\cos\theta)\dd\theta,
\qquad r\in\R.
\end{equation}
\begin{lemma}\label{lem:app-bessel-angular}
For every $r\in\R$, the Bessel function $J_0$ satisfies
\begin{equation}\label{eq:app-bessel-angular}
J_0(r)
=\frac1{2\pi}\int_0^{2\pi}e^{-ir\cos\theta}\dd\theta
=\frac1{2\pi}\int_0^{2\pi}\cos(r\cos\theta)\dd\theta.
\end{equation}
The same formula holds with $e^{ir\cos\theta}$ in place of $e^{-ir\cos\theta}$.
\end{lemma}

\begin{proof}
To prove both angular representations, it suffices by Euler's formula and \eqref{eq:bessel-definition} to show that $\int_0^{2\pi}\sin(r\cos\theta)\dd\theta=0$. Splitting this integral at $\pi$ and using $\cos(\theta+\pi)=-\cos\theta$, we obtain
\[
\int_0^{2\pi}\sin(r\cos\theta)\dd\theta
=\int_0^\pi\bigl[\sin(r\cos\theta)+\sin(r\cos(\theta+\pi))\bigr]\dd\theta=0.
\]
The exponential averages therefore equal $J_0(r)$ for either sign.
\end{proof}

We express the energy on a Fourier annulus as an integral of its radial density. We set $\mathbf e_\theta=(\cos\theta,\sin\theta)$. For a velocity field $u(t,x)$ and cutoffs $0<a<A<\infty$, polar coordinates give
\[
\frac1{2(2\pi)^2}\int_{a<|\xi|<A}|\wh u(t,\xi)|^2\dd\xi
=\frac1{2(2\pi)^2}\int_a^A k\int_0^{2\pi}|\wh u(t,k\mathbf e_\theta)|^2\dd\theta\dd k
=\int_a^A\sE(k,t)\dd k,
\]
where the radial kinetic-energy spectrum is defined by
\begin{equation}\label{eq:app-radial-spectrum}
\sE(k,t):=\frac{k}{2(2\pi)^2}
\int_0^{2\pi}|\wh u(t,k\mathbf e_\theta)|^2\dd\theta,
\qquad k>0.
\end{equation}

For the normalized three-vortex collapse \eqref{eq:configuration}--\eqref{eq:centers}, we define the time-independent spectral profile
\begin{equation}\label{eq:profiles}
F(q):=\sum_{k=1}^3\Gamma_k^2+
2\sum_{1\le i<j\le3}\Gamma_i\Gamma_jJ_0(qr_{ij}(0)),
\qquad q\ge0.
\end{equation}
The corresponding point-vortex spectrum is
\begin{equation}\label{eq:spectrum-profile}
\sE(k,t)=\frac1{4\pi k}F\!\left(k\sqrt{1-t}\right),
\qquad k>0,\quad 0\le t<1.
\end{equation}
This identity follows from Lemma~\ref{lem:app-spectrum} and the distance scaling $r_{ij}(t)=\sqrt{1-t}\,r_{ij}(0)$. In particular, the profile satisfies $F(k)=4\pi k\sE(k,0)$. The function $F$ gives the fixed spectral profile, and $\sE(k,t)$ gives the energy density in the physical wavenumber $k$. Time enters through the rescaled wavenumber $q=k\sqrt{1-t}$, for which $q\,r_{ij}(0)=k\,r_{ij}(t)$.

For a system of $n$ point vortices, we define the vorticity and its induced velocity by
\begin{equation}\label{eq:app-point-vortices}
\omega(t)=\sum_{i=1}^n\Gamma_i\delta_{z_i(t)},
\qquad u(t)=K*\omega(t),\qquad \wh\omega(t,\xi)=\sum_{i=1}^n\Gamma_i e^{-i\xi\cdot z_i(t)},
\end{equation}
where $K$ is the Biot--Savart kernel in \eqref{eq:NS}. The Fourier multiplier of $K$ gives
\begin{equation}\label{eq:BS-multiplier}
\wh u(t,\xi)=-i\frac{\xi^\perp}{|\xi|^2}\wh\omega(t,\xi),
\qquad \xi\ne0.
\end{equation}

\begin{lemma}
\label{lem:app-spectrum}
The radial kinetic-energy spectrum of the induced velocity is
\begin{equation}\label{eq:app-Novikov-spectrum}
\sE(k,t)
=\frac1{4\pi k}\left[
\sum_{i=1}^n\Gamma_i^2
+2\sum_{i<j}\Gamma_i\Gamma_j
J_0(kr_{ij}(t))\right],
\qquad k>0.
\end{equation}
\end{lemma}

\begin{proof}
To prove \eqref{eq:app-Novikov-spectrum}, we fix $k>0$ and a time at which the vortex positions are defined. By \eqref{eq:app-radial-spectrum}, it suffices to evaluate the angular integral of $|\wh u(t,k\mathbf e_\theta)|^2$. Combining \eqref{eq:BS-multiplier} with \eqref{eq:app-point-vortices}, we obtain
\[
\wh u(t,k\mathbf e_\theta)
=-\frac{i}{k}\mathbf e_\theta^\perp\wh\omega(t,k\mathbf e_\theta)
=-\frac{i}{k}\mathbf e_\theta^\perp\sum_{i=1}^n\Gamma_i e^{-ik\mathbf e_\theta\cdot z_i(t)}.
\]
Since $|\mathbf e_\theta^\perp|=1$, the squared modulus is $k^{-2}|\sum_i\Gamma_i e^{-ik\mathbf e_\theta\cdot z_i(t)}|^2$. Pairing the conjugate terms indexed by $(i,j)$ and $(j,i)$ gives
\[
|\wh u(t,k\mathbf e_\theta)|^2
=\frac1{k^2}\left[\sum_i\Gamma_i^2
+2\sum_{i<j}\Gamma_i\Gamma_j\cos\bigl(k\mathbf e_\theta\cdot(z_i(t)-z_j(t))\bigr)\right].
\]
The diagonal terms are independent of $\theta$, so it remains to identify the cosine integral for each pair. Writing $z_i(t)-z_j(t)=r_{ij}(t)\mathbf e_{\phi_{ij}}$ and setting $\vartheta=\theta-\phi_{ij}$, we use periodicity and \eqref{eq:bessel-definition} to obtain
\[
\int_0^{2\pi}\cos\bigl(k\mathbf e_\theta\cdot(z_i(t)-z_j(t))\bigr)\dd\theta
=\int_0^{2\pi}\cos\bigl(kr_{ij}(t)\cos\vartheta\bigr)\dd\vartheta
=2\pi J_0(kr_{ij}(t)).
\]
Each diagonal term also contributes a factor $2\pi$. Substituting these integrals into \eqref{eq:app-radial-spectrum} and using $\frac{k}{2(2\pi)^2}\frac{2\pi}{k^2}=\frac1{4\pi k}$ proves \eqref{eq:app-Novikov-spectrum}.
\end{proof}

\subsection{The mean-wavenumber prediction}
\label{sec:novikov-prediction}
\label{rem:app-mean-wavenumber}
We define the energy-weighted mean wavenumber by
\begin{equation}\label{eq:app-mean-wavenumber}
\ol k_{a,A}(t)
:=\frac{\displaystyle\int_a^A k\,\sE(k,t)\dd k}
{\displaystyle\int_a^A \sE(k,t)\dd k}.
\end{equation}
For the three-vortex collapse, the interaction moment is $\cM(t)=\frac1{2\pi}\sum_{i<j}\Gamma_i\Gamma_j/r_{ij}(t)$. The distance scaling \eqref{eq:distances} and Lemma~\ref{lem:app-interaction-moment} give $\cM(t)=\cM(0)/\sqrt{1-t}$, with $\cM(0)<0$ and $\cM'(t)<0$. We separate this time-dependent interaction from the individual-vortex contributions to describe the formal decrease of the mean wavenumber. With cutoff remainders $\cR_0$ and $\cR_1$ estimated in \eqref{eq:app-cutoff-remainders}, the exact formula is
\[
\ol k_{a,A}(t)
\overset{\text{Lemma~\ref{lem:app-spectral-moments}}}{=}
\frac{\displaystyle
\frac{A-a}{4\pi}\sum_{i=1}^3\Gamma_i^2
+\cM(t)+\cR_1(a,A,t)}
{\displaystyle
\frac{\log(A/a)}{4\pi}\sum_{i=1}^3\Gamma_i^2
+\cH+\cR_0(a,A,t)}.
\]
In the scale-separated regime $A^{-1}\ll r_{ij}(t)\ll a^{-1}$, omitting the cutoff remainders gives the approximation
\begin{equation}\label{eq:app-mean-wavenumber-collapse}
\ol k_{a,A}(t)\simeq
\frac{\displaystyle\frac{A-a}{4\pi}
\sum_{i=1}^3\Gamma_i^2+\cM(t)}
{\displaystyle\frac{\log(A/a)}{4\pi}
\sum_{i=1}^3\Gamma_i^2+\cH}.
\end{equation}
For fixed cutoffs with $A/a$ sufficiently large, Lemma~\ref{lem:app-interaction-moment} and the constant denominator imply that the approximation \eqref{eq:app-mean-wavenumber-collapse} decreases during collapse. This formal calculation gives the Novikov--Sedov prediction of a shift toward smaller wavenumbers and hence larger spatial scales while the cutoff separation holds.

\subsection{Finite-cutoff spectral moments}
\label{sec:novikov-sedov-details}

We justify the mean-wavenumber formula in Subsection~\ref{sec:novikov-prediction} by integrating the exact spectrum. We identify the finite interaction terms and estimate the errors at both cutoffs.

\subsubsection{The exact mean wavenumber}

\begin{lemma}\label{lem:app-mean-wavenumber-exact}
For a system of $n$ point vortices as in \eqref{eq:app-point-vortices}, we fix cutoffs $0<a<A<\infty$. When its denominator is positive, the mean wavenumber defined in \eqref{eq:app-mean-wavenumber} satisfies
\begin{equation}\label{eq:app-mean-wavenumber-exact}
\ol k_{a,A}(t)
=
\frac{\displaystyle
\frac{A-a}{4\pi}\sum_{i=1}^n\Gamma_i^2
+\frac1{2\pi}\sum_{i<j}\Gamma_i\Gamma_j
\int_a^A J_0(kr_{ij}(t))\dd k}
{\displaystyle
\frac{\log(A/a)}{4\pi}\sum_{i=1}^n\Gamma_i^2
+\frac1{2\pi}\sum_{i<j}\Gamma_i\Gamma_j
\int_a^A J_0(kr_{ij}(t))\frac{\dd k}{k}}.
\end{equation}
\end{lemma}

\begin{proof}
To prove \eqref{eq:app-mean-wavenumber-exact}, we fix $0<a<A<\infty$ and assume that the denominator in \eqref{eq:app-mean-wavenumber} is positive. It suffices to evaluate the two spectral integrals in that ratio. The formula \eqref{eq:app-Novikov-spectrum} reduces each evaluation to the individual-vortex term and a finite sum of pair terms. Integrating these terms over $[a,A]$ gives the identities
\begin{equation}\label{eq:app-finite-moments}
\begin{aligned}
\int_a^A k\sE(k,t)\dd k
&=\frac{A-a}{4\pi}\sum_{i=1}^n\Gamma_i^2
+\frac1{2\pi}\sum_{i<j}\Gamma_i\Gamma_j
\int_a^A J_0(kr_{ij}(t))\dd k,\\
\int_a^A\sE(k,t)\dd k
&=\frac{\log(A/a)}{4\pi}\sum_{i=1}^n\Gamma_i^2
+\frac1{2\pi}\sum_{i<j}\Gamma_i\Gamma_j
\int_a^A J_0(kr_{ij}(t))\frac{\dd k}{k}.
\end{aligned}
\end{equation}
The assumed positivity allows us to divide the first identity by the second. Substitution into \eqref{eq:app-mean-wavenumber} then proves \eqref{eq:app-mean-wavenumber-exact}.
\end{proof}

\subsubsection{Bessel integral estimates}

The leading terms in the finite-cutoff moments \eqref{eq:app-finite-moments} are determined by integrals of $J_0(kr)$ against $\dd k$ and $\dd k/k$. The following lemma evaluates these integrals and estimates the cutoff errors.

\begin{lemma}
\label{lem:app-bessel-integrals}
The improper integral of $J_0$ satisfies
\begin{equation}\label{eq:app-Bessel-mass}
\int_0^\infty J_0(s)\dd s=1.
\end{equation}
The constant $\kappa_0=\int_0^1\frac{J_0(s)-1}{s}\dd s+\int_1^\infty\frac{J_0(s)}s\dd s$ is finite. For every $r>0$ and $0<a<A<\infty$ with $ar\le1\le Ar$, we have the following estimates with universal implicit constants:
\begin{equation}\label{eq:app-Bessel-cutoffs}
\begin{split}
\left|\int_a^A J_0(kr)\dd k-\frac1r\right|
&\lesssim \left(a+\frac{(Ar)^{-1/2}}r\right),\\
\left|\int_a^A J_0(kr)\frac{\dd k}{k}
+\log(ar)-\kappa_0\right|
&\lesssim \left((ar)^2+(Ar)^{-3/2}\right).
\end{split}
\end{equation}
\end{lemma}

\begin{proof}
To prove \eqref{eq:app-Bessel-cutoffs}, we subtract the full Bessel integral for the first estimate and the constant $\kappa_0$ for the second, reducing both estimates to bounds on the omitted endpoint integrals. We first establish \eqref{eq:app-Bessel-mass}. To this end, we prove convergence of the improper integral, evaluate its exponentially damped version, and justify passage to the undamped limit.

For convergence and the upper-cutoff errors, we claim that, for $y\ge1$,
\begin{equation}\label{eq:app-Bessel-tails}
\int_y^\infty\frac{J_0(s)}{s^m}\dd s
=-\sqrt{\frac2\pi}\,y^{-m-1/2}\sin\left(y-\frac\pi4\right)
+O(y^{-m-3/2}),\qquad m=0,1.
\end{equation}
To prove this tail expansion, we use the standard Hankel expansion~\cite{DLMF},
\begin{equation}\label{eq:Hankel}
J_0(s)=\sqrt{\frac2{\pi s}}
\left[
\cos\left(s-\frac\pi4\right)
+\frac1{8s}\sin\left(s-\frac\pi4\right)
+O(s^{-2})
\right],\qquad s\to\infty.
\end{equation}
We fix $m\in\{0,1\}$. Integrating the leading cosine term in \eqref{eq:Hankel} by parts on $[y,R]$ and letting $R\to\infty$ gives
\[
\begin{aligned}
\int_y^\infty s^{-m-1/2}\cos\left(s-\frac\pi4\right)\dd s
&=-y^{-m-1/2}\sin\left(y-\frac\pi4\right)\\
&\quad+\left(m+\frac12\right)\int_y^\infty s^{-m-3/2}\sin\left(s-\frac\pi4\right)\dd s.
\end{aligned}
\] The remainder in \eqref{eq:app-Bessel-tails} is therefore bounded by a constant multiple of $(\mathrm I)+(\mathrm{II})$, where $(\mathrm I)=|\int_y^\infty s^{-m-3/2}\sin(s-\pi/4)\dd s|$ and $(\mathrm{II})=\int_y^\infty s^{-m-5/2}\dd s$. A second integration by parts yields
\[
\begin{aligned}
\int_y^\infty s^{-m-3/2}\sin\left(s-\frac\pi4\right)\dd s
&=y^{-m-3/2}\cos\left(y-\frac\pi4\right)\\
&\quad-\left(m+\frac32\right)\int_y^\infty s^{-m-5/2}\cos\left(s-\frac\pi4\right)\dd s.
\end{aligned}
\] Using $|\cos|\le1$, we obtain $(\mathrm I)\le y^{-m-3/2}+(m+\frac32)\int_y^\infty s^{-m-5/2}\dd s=2y^{-m-3/2}$ and $(\mathrm{II})=y^{-m-3/2}/(m+\frac32)\lesssim y^{-m-3/2}$. These bounds prove \eqref{eq:app-Bessel-tails}; its case $m=0$ establishes convergence of the Bessel integral, and its case $m=1$ will control the logarithmic upper-cutoff error.

To evaluate this integral, we fix $\epsilon>0$ and evaluate the damped integral. Absolute convergence justifies Fubini's theorem, so \eqref{eq:bessel-definition}, symmetry, and the substitution $u=\tan\theta$ give
\[
\begin{aligned}
\int_0^\infty e^{-\epsilon s}J_0(s)\dd s
&=\frac1{2\pi}\int_0^{2\pi}\left(\int_0^\infty e^{-\epsilon s}\cos(s\cos\theta)\dd s\right)\dd\theta
=\frac1{2\pi}\int_0^{2\pi}\frac{\epsilon}{\epsilon^2+\cos^2\theta}\dd\theta\\
&=\frac{2\epsilon}{\pi}\int_0^{\pi/2}\frac{\dd\theta}{\epsilon^2+\cos^2\theta}
=\frac{2\epsilon}{\pi}\int_0^\infty\frac{\dd u}{1+\epsilon^2+\epsilon^2u^2}
=\frac1{\sqrt{1+\epsilon^2}}.
\end{aligned}
\] The last integral is evaluated using the antiderivative $[\epsilon\sqrt{1+\epsilon^2}]^{-1}\arctan(\epsilon u/\sqrt{1+\epsilon^2})$. To remove the damping, we set $B(x)=\int_0^xJ_0(s)\dd s$. The convergence already proved makes $B$ bounded and gives it a finite limit at infinity. Integration by parts and the substitution $v=\epsilon s$ yield
\[
\int_0^\infty e^{-\epsilon s}J_0(s)\dd s
=\epsilon\int_0^\infty e^{-\epsilon s}B(s)\dd s
=\int_0^\infty e^{-v}B(v/\epsilon)\dd v.
\]
The last integrand is dominated by a constant multiple of $e^{-v}$. By dominated convergence, the improper integral equals $\lim_{\epsilon\to0^+}(1+\epsilon^2)^{-1/2}=1$, proving \eqref{eq:app-Bessel-mass}.

To prove the first cutoff estimate, we fix $r>0$ and $0<a<A$ with $ar\le1\le Ar$. Subtracting \eqref{eq:app-Bessel-mass} and setting $s=kr$ gives
\[
\int_a^A J_0(kr)\dd k-\frac1r
=-\underbrace{\frac1r\int_0^{ar}J_0(s)\dd s}_{(\mathrm I)}
-\underbrace{\frac1r\int_{Ar}^\infty J_0(s)\dd s}_{(\mathrm{II})}.
\]
The angular definition gives $|J_0|\le1$, hence $|(\mathrm I)|\le r^{-1}\int_0^{ar}|J_0(s)|\dd s\le r^{-1}\int_0^{ar}1\dd s=a$. Since $Ar\ge1$, the case $m=0$ of \eqref{eq:app-Bessel-tails} gives $|(\mathrm{II})|\lesssim r^{-1}(Ar)^{-1/2}$. Combining these bounds with the endpoint decomposition proves the first estimate in \eqref{eq:app-Bessel-cutoffs}.

For the logarithmic estimate, we first justify subtraction of $\kappa_0$. The expansion $J_0(s)-1=O(s^2)$ at zero and \eqref{eq:Hankel} at infinity give $\int_0^1|J_0(s)-1|s^{-1}\dd s\lesssim\int_0^1s\dd s\lesssim1$ and $\int_1^\infty|J_0(s)|s^{-1}\dd s\lesssim\int_1^\infty s^{-3/2}\dd s\lesssim1$. Both defining integrals are thus absolutely convergent, with $|\kappa_0|\lesssim1$. Setting $s=kr$ and subtracting its definition, we obtain
\[
\int_a^A J_0(kr)\frac{\dd k}{k}+\log(ar)-\kappa_0
=-\underbrace{\int_0^{ar}\frac{J_0(s)-1}{s}\dd s}_{(\mathrm I)}
-\underbrace{\int_{Ar}^\infty\frac{J_0(s)}s\dd s}_{(\mathrm{II})}.
\]
To bound the lower endpoint, we integrate the expansion $J_0(s)-1=-s^2/4+O(s^4)$ and obtain
\begin{equation}\label{eq:app-Bessel-origin}
\int_0^x\frac{J_0(s)-1}{s}\dd s=-\frac{x^2}{8}+O(x^4),\qquad 0<x\le1.
\end{equation}
Since $ar\le1\le Ar$, this expansion gives $|(\mathrm I)|\lesssim(ar)^2$, while the case $m=1$ of \eqref{eq:app-Bessel-tails} gives $|(\mathrm{II})|\lesssim(Ar)^{-3/2}$. Substitution into the endpoint decomposition proves the second estimate in \eqref{eq:app-Bessel-cutoffs}. All constants are independent of $a$, $A$, and $r$.
\end{proof}

\subsubsection{Spectral moments and the mean wavenumber}

Subtracting the explicit leading terms in \eqref{eq:app-finite-moments} separates the time-dependent interaction from the individual-vortex contributions. The following lemma identifies the finite interaction terms and bounds the cutoff errors.

\begin{lemma}
\label{lem:app-spectral-moments}
For the normalized three-vortex configuration \eqref{eq:configuration}, we denote by $\cM(t)$ the interaction moment defined in Subsection~\ref{sec:novikov-prediction}. The Hamiltonian $\cH$ is defined in \eqref{eq:transfer-constant}. If $0<a<A<\infty$ and $a\,r_{ij}(t)\le1\le A\,r_{ij}(t)$ for every $i<j$, then the spectral moments satisfy
\begin{equation}\label{eq:app-spectral-moments}
\begin{split}
\int_a^A k\sE(k,t)\dd k
&=\frac{A-a}{4\pi}\sum_{i=1}^3\Gamma_i^2
+\cM(t)+\cR_1(a,A,t),\\
\int_a^A\sE(k,t)\dd k
&=\frac{\log(A/a)}{4\pi}\sum_{i=1}^3\Gamma_i^2
+\cH+\cR_0(a,A,t),
\end{split}
\end{equation}
where the remainders satisfy
\begin{equation}\label{eq:app-cutoff-remainders}
\begin{split}
|\cR_1(a,A,t)|
&\lesssim \sum_{i<j}|\Gamma_i\Gamma_j|
\left(a+
\frac{(A r_{ij}(t))^{-1/2}}{r_{ij}(t)}\right),\\
|\cR_0(a,A,t)|
&\lesssim \sum_{i<j}|\Gamma_i\Gamma_j|
\left((a r_{ij}(t))^2
+(A r_{ij}(t))^{-3/2}\right).
\end{split}
\end{equation}
Both remainders therefore tend to zero as $a\to0$ and $A\to\infty$ for each fixed $t<1$. The renormalized first spectral moment satisfies
\begin{equation}\label{eq:app-spectral-moment}
\cM(t)
=
\int_0^\infty k\left[
\sE(k,t)
-\frac1{4\pi k}\sum_{i=1}^3\Gamma_i^2
\right]\dd k,
\end{equation}
where the integral is improper.
\end{lemma}

\begin{proof}
To prove \eqref{eq:app-spectral-moments} and \eqref{eq:app-cutoff-remainders}, we fix $0\le t<1$ and cutoffs $0<a<A$ satisfying the hypotheses. We define the remainders by subtracting the asserted leading terms from \eqref{eq:app-finite-moments}; it then suffices to bound these differences by the stated endpoint errors. The improper-integral formula \eqref{eq:app-spectral-moment} will follow from the first remainder bound as $a\to0$ and $A\to\infty$.

For the first moment, the definition of $\cM(t)$ and the Bessel integral identity give
\[
\begin{aligned}
\cR_1(a,A,t)
&:=\int_a^A k\sE(k,t)\dd k-\frac{A-a}{4\pi}\sum_i\Gamma_i^2-\cM(t)\\
&=-\underbrace{\frac1{2\pi}\sum_{i<j}\frac{\Gamma_i\Gamma_j}{r_{ij}(t)}\int_0^{ar_{ij}(t)}J_0(s)\dd s}_{(\mathrm I)}
-\underbrace{\frac1{2\pi}\sum_{i<j}\frac{\Gamma_i\Gamma_j}{r_{ij}(t)}\int_{Ar_{ij}(t)}^\infty J_0(s)\dd s}_{(\mathrm{II})}.
\end{aligned}
\] The bound $|J_0|\le1$ gives $|(\mathrm I)|\le\frac1{2\pi}\sum_{i<j}\frac{|\Gamma_i\Gamma_j|}{r_{ij}(t)}\int_0^{ar_{ij}(t)}1\dd s=\frac a{2\pi}\sum_{i<j}|\Gamma_i\Gamma_j|$. Since $Ar_{ij}(t)\ge1$ for every pair, \eqref{eq:app-Bessel-tails} gives
\[
|(\mathrm{II})|
\le\frac1{2\pi}\sum_{i<j}\frac{|\Gamma_i\Gamma_j|}{r_{ij}(t)}\left|\int_{Ar_{ij}(t)}^\infty J_0(s)\dd s\right|
\lesssim\sum_{i<j}|\Gamma_i\Gamma_j|\frac{(Ar_{ij}(t))^{-1/2}}{r_{ij}(t)}.
\]
The inequality $|\cR_1|\le|(\mathrm I)|+|(\mathrm{II})|$ now proves the first bound in \eqref{eq:app-cutoff-remainders} and the corresponding moment expansion.

For the energy integral, the remainder is
\[
\begin{aligned}
\cR_0(a,A,t)
&:=\int_a^A\sE(k,t)\dd k-\frac{\log(A/a)}{4\pi}\sum_i\Gamma_i^2-\cH\\
&=\frac1{2\pi}\sum_{i<j}\Gamma_i\Gamma_j\int_a^A J_0(kr_{ij}(t))\frac{\dd k}{k}-\cH.
\end{aligned}
\] To apply \eqref{eq:app-Bessel-cutoffs}, we add $\log(a r_{ij}(t))-\kappa_0$ inside each pair integral and check that these corrections cancel after summation. The distance law gives $\log(a r_{ij}(t))=\log a+\frac12\log(1-t)+\log r_{ij}(0)$, so
\[
\begin{aligned}
\cR_0(a,A,t)
&=\frac1{2\pi}\sum_{i<j}\Gamma_i\Gamma_j
\left[\int_a^A J_0(kr_{ij}(t))\frac{\dd k}{k}+\log(a r_{ij}(t))-\kappa_0\right]\\
&\quad+\frac{\kappa_0-\log a-\frac12\log(1-t)}{2\pi}\sum_{i<j}\Gamma_i\Gamma_j
-\left[\frac1{2\pi}\sum_{i<j}\Gamma_i\Gamma_j\log r_{ij}(0)+\cH\right].
\end{aligned}
\] The circulation identity \eqref{eq:invariants} cancels the second term, and \eqref{eq:transfer-constant} cancels the third. Since $ar_{ij}(t)\le1\le Ar_{ij}(t)$ for every pair, applying \eqref{eq:app-Bessel-cutoffs} to the remaining sum yields
\[
\begin{aligned}
|\cR_0(a,A,t)|
&\le\frac1{2\pi}\sum_{i<j}|\Gamma_i\Gamma_j|
\left|\int_a^A J_0(kr_{ij}(t))\frac{\dd k}{k}+\log(a r_{ij}(t))-\kappa_0\right|\\
&\lesssim\sum_{i<j}|\Gamma_i\Gamma_j|\left((ar_{ij}(t))^2+(Ar_{ij}(t))^{-3/2}\right).
\end{aligned}
\] This proves the energy expansion and its error bound. At the fixed time under consideration, all pair distances are positive constants, so both remainders tend to zero as $a\to0$ and $A\to\infty$.

To identify the renormalized improper integral, we subtract the individual-vortex contribution from the first moment expansion and obtain
\[
\int_a^A k\left[\sE(k,t)-\frac1{4\pi k}\sum_i\Gamma_i^2\right]\dd k
=\cM(t)+\cR_1(a,A,t).
\]
Since $\cR_1$ tends to zero, passage to the limit at both endpoints proves \eqref{eq:app-spectral-moment}.
\end{proof}

\begin{lemma}
\label{lem:app-interaction-moment}
For the normalized collapsing configuration \eqref{eq:configuration}, the initial interaction moment satisfies $\cM(0)<0$. Its evolution is given by
\begin{equation}\label{eq:app-contracting-moment}
\cM(t)=\frac{\cM(0)}{\sqrt{1-t}},
\qquad
\cM'(t)
=\frac{\cM(0)}{2(1-t)^{3/2}}
=\frac{\cM(t)}{2(1-t)}<0,
\qquad 0\le t<1.
\end{equation}
\end{lemma}

\begin{proof}
To prove \eqref{eq:app-contracting-moment}, it suffices to establish $\cM(0)<0$. Indeed, substituting \eqref{eq:distances} into the definition gives $\cM(t)=\frac1{2\pi}\sum_{i<j}\Gamma_i\Gamma_j/(\sqrt{1-t}\,r_{ij}(0))=\cM(0)/\sqrt{1-t}$. Differentiation yields $\cM'(t)=\cM(0)/[2(1-t)^{3/2}]=\cM(t)/[2(1-t)]$, so both strict signs follow from that of $\cM(0)$.

To determine this sign, we use the labeling of Lemma~\ref{lem:distance-geometry}, for which $\Gamma_1\Gamma_2>0$. The circulation identity gives $\Gamma_3=-\Gamma_1\Gamma_2/(\Gamma_1+\Gamma_2)$, and substitution into the definition of $\cM(0)$ yields
\[
\frac{2\pi\cM(0)}{\Gamma_1\Gamma_2}
=\frac1{r_{12}(0)}
-\frac{\Gamma_1}{\Gamma_1+\Gamma_2}\frac1{r_{13}(0)}
-\frac{\Gamma_2}{\Gamma_1+\Gamma_2}\frac1{r_{23}(0)}.
\]
Since $\Gamma_1\Gamma_2>0$, the desired negativity is equivalent to $1/r_{12}(0)<\frac{\Gamma_1}{\Gamma_1+\Gamma_2}/r_{13}(0)+\frac{\Gamma_2}{\Gamma_1+\Gamma_2}/r_{23}(0)$. To establish this inequality, we use the strict convexity of $s\mapsto s^{-1/2}$ and the squared-distance relation \eqref{eq:second-moment}, which has the form
\[
r_{12}(0)^2
=\frac{\Gamma_1}{\Gamma_1+\Gamma_2}r_{13}(0)^2
+\frac{\Gamma_2}{\Gamma_1+\Gamma_2}r_{23}(0)^2.
\]
The two coefficients belong to $(0,1)$ and sum to one, while Lemma~\ref{lem:distance-geometry} gives $r_{13}(0)\ne r_{23}(0)$. Strict convexity therefore implies
\[
\begin{aligned}
\frac1{r_{12}(0)}
&=\left[\frac{\Gamma_1}{\Gamma_1+\Gamma_2}r_{13}(0)^2
+\frac{\Gamma_2}{\Gamma_1+\Gamma_2}r_{23}(0)^2\right]^{-1/2}\\
&<\frac{\Gamma_1}{\Gamma_1+\Gamma_2}\frac1{r_{13}(0)}
+\frac{\Gamma_2}{\Gamma_1+\Gamma_2}\frac1{r_{23}(0)}.
\end{aligned}
\] This proves $\cM(0)<0$, and the scaling and derivative identities give all remaining conclusions.
\end{proof}

\section{Hamiltonian and spectral profile}\label{sec:spectrum}

We prove that the transfer constant $\cH$ is negative and establish the low-frequency
bound needed to estimate the Gaussian energy. The circulation and distance identities below
yield the cancellation of the quadratic term.

\begin{proposition}\label{thm:spectral-positivity}
For the configuration \eqref{eq:configuration}, the Hamiltonian
satisfies $\cH<0$. The value of $\cH$ is invariant
under multiplication of all distances $r_{ij}(0)$ by a common
positive factor.
\end{proposition}

To prove $\cH<0$, we rewrite its logarithmic pair sum in terms of a convex combination
of two distinct squared distances and apply strict concavity. The following identities determine
this combination.

\begin{lemma}[{\cite{NS79}}]\label{lem:invariants}
The configuration \eqref{eq:configuration} satisfies
\begin{equation}\label{eq:invariants}
\sum_{i<j}\Gamma_i\Gamma_j=0,
\qquad \sum_i\Gamma_i z_i^0=0,
\qquad \sum_i\Gamma_i|z_i^0|^2=0,
\end{equation}
and
\begin{equation}\label{eq:second-moment}
\sum_{i<j}\Gamma_i\Gamma_jr_{ij}(0)^2=0.
\end{equation}
The circulations have both signs, and the total circulation satisfies
\begin{equation}\label{eq:circulation-square}
\Gtot^2=\sum_i\Gamma_i^2>0.
\end{equation}
\end{lemma}

For completeness, we give a proof.
\begin{proof}
To prove the distance identity \eqref{eq:second-moment} and the
circulation formula \eqref{eq:circulation-square}, it suffices to
establish the three cancellations in \eqref{eq:invariants}.
Indeed, these cancellations give
\[
\begin{aligned}
\sum_{i<j}\Gamma_i\Gamma_jr_{ij}(0)^2
&=\Gtot\sum_i\Gamma_i|z_i^0|^2-\left|\sum_i\Gamma_i z_i^0\right|^2=0,\\
\Gtot^2
&=\sum_i\Gamma_i^2+2\sum_{i<j}\Gamma_i\Gamma_j=\sum_i\Gamma_i^2>0.
\end{aligned}
\]
If the three nonzero circulations had the same sign, then $\sum_{i<j}\Gamma_i\Gamma_j>0$, contradicting \eqref{eq:invariants}.

To prove the pair-circulation cancellation, it is enough to show
that the quadratic moment vanishes. Taking the scalar product of
\eqref{eq:configuration} with $\Gamma_i(z_i^0)^\perp$, summing over
$i$, and pairing the two orders of each vortex pair, we obtain
\[
\begin{aligned}
0
&=-\beta\sum_i\Gamma_i|z_i^0|^2
+\sum_{i<j}\Gamma_i\Gamma_j
(z_i^0-z_j^0)^\perp\cdot K(z_i^0-z_j^0)\\
&=-\beta\sum_i\Gamma_i|z_i^0|^2
+\frac1{2\pi}\sum_{i<j}\Gamma_i\Gamma_j.
\end{aligned}
\]
To show that the quadratic moment vanishes, we take the scalar
product of \eqref{eq:configuration} with $\Gamma_i z_i^0$.
The same pairing and the identity $x\cdot K(x)=0$ give
\[
0=\frac12\sum_i\Gamma_i|z_i^0|^2+\sum_{i<j}\Gamma_i\Gamma_j(z_i^0-z_j^0)\cdot K(z_i^0-z_j^0)=\frac12\sum_i\Gamma_i|z_i^0|^2.
\]
Substitution into the preceding identity now proves
$\sum_{i<j}\Gamma_i\Gamma_j=0$.

It remains to prove the vector cancellation in
\eqref{eq:invariants}. Multiplication of \eqref{eq:configuration}
by $\Gamma_i$ and summation over $i$ give $\frac12\sum_i\Gamma_i z_i^0-\beta(\sum_i\Gamma_i z_i^0)^\perp=0$, because the oddness of $K$ cancels the two contributions of each pair. Taking the scalar product with $\sum_i\Gamma_i z_i^0$ yields $\frac12|\sum_i\Gamma_i z_i^0|^2=0$. This proves all three cancellations in \eqref{eq:invariants}.
The two expansions at the beginning of the proof then yield
\eqref{eq:second-moment} and
$\Gtot^2=\sum_i\Gamma_i^2>0$, as required.
\end{proof}

\begin{lemma}\label{lem:distance-geometry}
The distances $r_{12}(0),r_{13}(0),r_{23}(0)$ are pairwise
distinct. We relabel the vortices so that $\Gamma_1\Gamma_2>0$ and set $\theta=\Gamma_1/(\Gamma_1+\Gamma_2)\in(0,1)$, $X=r_{12}(0)^2$, $Y=r_{13}(0)^2$, and $Z=r_{23}(0)^2$. These quantities satisfy
\begin{equation}\label{eq:distance-convexity}
X=\theta Y+(1-\theta)Z,\qquad Y\ne Z.
\end{equation}
\end{lemma}

\begin{proof}
To prove both assertions, it suffices to establish the convex
combination in \eqref{eq:distance-convexity} and exclude $Y=Z$.
Indeed, a strict convex combination of distinct $Y$ and $Z$
lies strictly between them, so the three squared distances,
and hence the three distances, are pairwise distinct.

We use the relabeling in the statement, which is possible because
both signs occur among the three nonzero circulations.
To derive the convex combination, we solve the pair-circulation
identity in \eqref{eq:invariants} for $\Gamma_3$, obtaining $\Gamma_3=-\Gamma_1\Gamma_2/(\Gamma_1+\Gamma_2)$. The denominator is nonzero because $\Gamma_1$ and $\Gamma_2$ have the same sign. Substitution into \eqref{eq:second-moment}
then gives
\[
0=\Gamma_1\Gamma_2X+\Gamma_1\Gamma_3Y+\Gamma_2\Gamma_3Z
=\Gamma_1\Gamma_2\left[X-\frac{\Gamma_1Y+\Gamma_2Z}{\Gamma_1+\Gamma_2}\right].
\]
Division by $\Gamma_1\Gamma_2\ne0$ proves
$X=\theta Y+(1-\theta)Z$, with $0<\theta<1$.

To exclude equality of the remaining two squared distances,
we assume that $Y=Z$. The convex combination then implies
$X=Y=Z=r^2$ for some $r>0$. The collapse law
\eqref{eq:distances} therefore requires
\[
\left.\frac{\dd}{\dd t}|z_i(t)-z_j(t)|^2\right|_{t=0}
=-r^2<0
\qquad(i\ne j).
\]
We obtain a contradiction by showing that the point-vortex
velocities make this derivative vanish. At the equilateral
configuration, their equations give
\[
\begin{aligned}
\dot z_i(0)
&=\frac1{2\pi r^2}\sum_{k\ne i}\Gamma_k(z_i^0-z_k^0)^\perp\\
&=\frac1{2\pi r^2}
\left(\Gtot z_i^0-\sum_k\Gamma_kz_k^0\right)^\perp.
\end{aligned}
\]
Subtracting the velocities of two vortices gives $\dot z_i(0)-\dot z_j(0)=\frac{\Gtot}{2\pi r^2}(z_i^0-z_j^0)^\perp$. Consequently, differentiation of the squared distance yields
\[
\left.\frac{\dd}{\dd t}|z_i(t)-z_j(t)|^2\right|_{t=0}
=\frac{\Gtot}{\pi r^2}
(z_i^0-z_j^0)\cdot(z_i^0-z_j^0)^\perp=0,
\]
contradicting the collapse law. Thus $Y\ne Z$, and the strict
inequalities $\min\{Y,Z\}<X<\max\{Y,Z\}$ prove the claim.
\end{proof}

\begin{proof}[Proof of Proposition~\ref{thm:spectral-positivity}]
To prove $\cH<0$, it suffices to express $-4\pi\cH$
as a positive multiple of the strict concavity gap for the logarithm.
We use the labeling of Lemma~\ref{lem:distance-geometry} and set $p=\Gamma_1\Gamma_2>0$, $\theta=\Gamma_1/(\Gamma_1+\Gamma_2)\in(0,1)$, $Y=r_{13}(0)^2$, and $Z=r_{23}(0)^2$. Since that lemma gives $Y\ne Z$, strict concavity will prove
the sign assertion once we establish the representation
\begin{equation}\label{eq:Jensen}
\cH=-\frac{p}{4\pi}
\log\frac{\theta Y+(1-\theta)Z}{Y^\theta Z^{1-\theta}}<0.
\end{equation}
To derive the equality in \eqref{eq:Jensen}, we first express
the pair products in terms of $p$ and $\theta$. The circulation
identity in \eqref{eq:invariants} gives $\Gamma_3=-p/(\Gamma_1+\Gamma_2)$, $\Gamma_1\Gamma_3=-p\theta$, and $\Gamma_2\Gamma_3=-p(1-\theta)$. Substitution into \eqref{eq:transfer-constant}, followed by
\eqref{eq:distance-convexity}, yields
\[
\begin{aligned}
-4\pi\cH
&=p\left[\log r_{12}(0)^2-\theta\log Y-(1-\theta)\log Z\right]\\
&=p\left[\log\bigl(\theta Y+(1-\theta)Z\bigr)
-\theta\log Y-(1-\theta)\log Z\right].
\end{aligned}
\]
This is the required representation, and strict concavity gives
the inequality in \eqref{eq:Jensen}.

To prove dilation invariance, it suffices to show that replacing
each distance by $s$ times that distance adds zero to the
Hamiltonian. We fix $s>0$ and use \eqref{eq:invariants} to compute
\[
-\frac1{2\pi}\sum_{i<j}\Gamma_i\Gamma_j\log(sr_{ij}(0))
=\cH-\frac{\log s}{2\pi}\sum_{i<j}\Gamma_i\Gamma_j=\cH.
\]
This proves the second assertion.
\end{proof}

The definition \eqref{eq:profiles} and the angular identity
\eqref{eq:app-bessel-angular} yield, for every $q\ge0$,
\begin{equation}\label{eq:F-angular}
F(q)=\frac1{2\pi}\int_0^{2\pi}
\left|\sum_j\Gamma_j e^{-iq \mathbf e_\theta\cdot z_j^0}\right|^2\dd\theta.
\end{equation}
In particular, the profile satisfies $F(0)=(\sum_j\Gamma_j)^2=\Gtot^2$.
We use the following quartic bound in the interaction estimates.

\begin{lemma}\label{lem:profile-origin}
For every $q\ge0$, the spectral profile satisfies
\begin{equation}\label{eq:F-origin}
\left|F(q)-\Gtot^2\right|
\le\frac{q^4}{32}
\sum_{i<j}|\Gamma_i\Gamma_j|r_{ij}(0)^4.
\end{equation}
\end{lemma}

\begin{proof}
To prove \eqref{eq:F-origin} for a fixed $q\ge0$, we reduce its
left-hand side to the fourth-order remainder of $J_0$.
The definition \eqref{eq:profiles}, the identity
$F(0)=\Gtot^2$, and the quadratic cancellation
\eqref{eq:second-moment} give
\[
F(q)-\Gtot^2=2\sum_{i<j}\Gamma_i\Gamma_j\bigl[J_0(qr_{ij}(0))-1\bigr]=2\sum_{i<j}\Gamma_i\Gamma_j\left[J_0(qr_{ij}(0))-1+\frac{q^2r_{ij}(0)^2}{4}\right].
\]
The triangle inequality therefore reduces the desired bound to the scalar estimate $|J_0(z)-1+z^2/4|\le |z|^4/64$ for $z\in\R$. To prove it, we use the angular representation
\eqref{eq:bessel-definition} to reduce the Bessel remainder
to the cosine remainder. Since the angular average of
$\cos^2\theta$ is $1/2$, we have
\[
J_0(z)-1+\frac{z^2}{4}
=\frac1{2\pi}\int_0^{2\pi}
\left[\cos(z\cos\theta)-1+\frac{z^2\cos^2\theta}{2}\right]
\dd\theta.
\]
Taylor's theorem gives $|\cos w-1+w^2/2|\le |w|^4/24$ for $w\in\R$. Applying this bound with $w=z\cos\theta$ and using the angular average of $\cos^4\theta$, we obtain
\[
\left|J_0(z)-1+\frac{z^2}{4}\right|
\le\frac{|z|^4}{24}\,\frac1{2\pi}\int_0^{2\pi}\cos^4\theta\dd\theta
=\frac{|z|^4}{24}\cdot\frac38=\frac{|z|^4}{64}.
\]
The scalar estimate applied to each pair in the reduced profile
identity now yields
\[
\left|F(q)-\Gtot^2\right|
\le2\sum_{i<j}|\Gamma_i\Gamma_j|\frac{q^4r_{ij}(0)^4}{64}
=\frac{q^4}{32}\sum_{i<j}|\Gamma_i\Gamma_j|r_{ij}(0)^4,
\]
which proves \eqref{eq:F-origin}.
\end{proof}

\section{Energy identities and viscous dissipation}
\label{sec:estimates}

We derive the energy identity above a fixed Fourier cutoff and bound the viscous loss on the
shrinking observation interval.

The vorticity conserves total circulation and satisfies
\begin{equation}\label{eq:L1-bound}
\int_{\R^2}\omega^\nu(t,x)\dd x=\Gtot,
\qquad \|\omega^\nu(t)\|_{L^1(\R^2)}\le\Gabs,
\qquad t>0.
\end{equation}
The identity \eqref{eq:circulation-square} implies
$\Gtot\ne0$, so the total kinetic energy diverges
at zero frequency. The energy above any positive cutoff is finite: for every $N,t>0$,
\[
\int_{\R^2}\ind_{\{|\xi|>N\}}|\wh{u^\nu}(t)|^2\dd\xi
\le\frac1{N^2}\int_{\R^2}
|\wh{\omega^\nu}(t)|^2\dd\xi<\infty.
\]

\begin{lemma}\label{lem:balances}
For every $\nu,N,t>0$, the energy above the cutoff satisfies
\begin{equation}\label{eq:high-balance}
\begin{aligned}
&\frac1{2(2\pi)^2}\frac{\dd}{\dd t}
\int_{\R^2}\ind_{\{|\xi|>N\}}
|\wh u^\nu(t)|^2\dd\xi
+\frac{\nu}{(2\pi)^2}\int_{\R^2}
\ind_{\{|\xi|>N\}}|\wh\omega^\nu(t)|^2\dd\xi\\
&\qquad=-\frac1{(2\pi)^2}\Real
\int_{\R^2}\ind_{\{|\xi|>N\}}
\ol{\wh u^\nu(t)}
\cdot\wh{(u^\nu\cdot\nabla)u^\nu}(t)\dd\xi.
\end{aligned}
\end{equation}
Here $|\wh\omega^\nu(t,\xi)|^2=|\xi|^2|\wh u^\nu(t,\xi)|^2$
by \eqref{eq:BS-multiplier}.
\end{lemma}

\begin{proof}
To prove \eqref{eq:high-balance} for fixed $\nu,N>0$, it suffices
to establish the Fourier energy identity
\[
\frac12\partial_t|\wh u^\nu|^2
+\nu|\xi|^2|\wh u^\nu|^2
=-\Real\bigl[
\ol{\wh u^\nu}\cdot
\wh{(u^\nu\cdot\nabla)u^\nu}\bigr],
\]
for almost every nonzero frequency at positive times. Indeed,
multiplication by $(2\pi)^{-2}\ind_{\{|\xi|>N\}}$ and
integration in $\xi$ give the claimed balance, because
\eqref{eq:BS-multiplier} identifies
$|\xi|^2|\wh u^\nu|^2$ with $|\wh\omega^\nu|^2$.

To establish the Fourier identity, we take the Fourier transform
of the velocity equation associated with \eqref{eq:NS}.
Writing $p^\nu$ for the pressure, we obtain
\[
\partial_t\wh u^\nu
+\wh{(u^\nu\cdot\nabla)u^\nu}
+i\xi\wh p^\nu
=-\nu|\xi|^2\wh u^\nu,
\qquad \xi\cdot\wh u^\nu=0.
\]
Taking the scalar product with $\ol{\wh u^\nu}$ and then the real part gives the required Fourier identity: the pressure contribution is $-\Real[i\wh p^\nu\,\xi\cdot\ol{\wh u^\nu}]=0$ because $\xi\cdot\ol{\wh u^\nu}=0$.

The operations above are justified on every compact time
interval in $(0,\infty)$ by the positive-time bounds
$u^\nu\in L^\infty(\R^2)$ and
$\omega^\nu\in H^1(\R^2)$, together with the cutoff
away from $\xi=0$. In particular, the Biot--Savart multiplier
gives $\nabla u^\nu\in L^2(\R^2)$, so
$(u^\nu\cdot\nabla)u^\nu\in L^2(\R^2)$.
Thus integration of the Fourier identity yields
\eqref{eq:high-balance}.
\end{proof}

The right-hand side of \eqref{eq:high-balance} gives the rate of nonlinear energy transfer into
$\{|\xi|>N\}$; a negative value indicates transfer toward lower frequencies. The nonnegative
viscous term gives the rate of energy loss by diffusion.

The vorticity smoothing estimate of Gallay and Wayne gives the following bound.

\begin{lemma}[{\cite{GW05}}]\label{lem:dissipation-bound}
For every $\nu,t>0$, the vorticity satisfies
\begin{equation}\label{eq:dissipation-bound}
\nu\|\omega^\nu(t)\|_{L^2(\R^2)}^2
\lesssim\frac{\Gabs^2}{t}.
\end{equation}
\end{lemma}

Integration of \eqref{eq:dissipation-bound} in time yields the following corollary.

\begin{corollary}\label{cor:vanishing-viscous-dissipation}
For every $\nu>0$ and $0<a<b$, the integrated viscous dissipation satisfies
\begin{equation}\label{eq:integrated-dissipation-bound}
0\le\nu\int_a^b\|\omega^\nu(t)\|_{L^2(\R^2)}^2\dd t
\lesssim\Gabs^2\log\frac ba.
\end{equation}
In particular, if $0<a_\nu<b_\nu$ and $b_\nu/a_\nu\to1$ as
$\nu\to0$, then
\begin{equation}\label{eq:vanishing-viscous-dissipation}
\nu\int_{a_\nu}^{b_\nu}
\|\omega^\nu(t)\|_{L^2(\R^2)}^2\dd t
\longrightarrow0.
\end{equation}
By Plancherel's identity, both conclusions remain valid uniformly for $N>0$ when the squared $L^2$
norm is replaced by $(2\pi)^{-2}\int_{|\xi|>N}|\wh\omega^\nu(t,\xi)|^2\dd\xi$.
\end{corollary}

\begin{proof}
To prove \eqref{eq:integrated-dissipation-bound}, it suffices to
integrate the smoothing estimate \eqref{eq:dissipation-bound}
over the interval $[a,b]$. Applying that estimate, we obtain
\[
0\le\nu\int_a^b\|\omega^\nu(t)\|_{L^2(\R^2)}^2\dd t
\lesssim\Gabs^2\int_a^b\frac{\dd t}{t}
=\Gabs^2\log\frac ba.
\]
The vanishing assertion \eqref{eq:vanishing-viscous-dissipation}
therefore follows from $\log(b_\nu/a_\nu)\to0$.

To obtain the corresponding conclusions above a Fourier cutoff,
it suffices to bound the restricted Fourier integral by the full
$L^2$ norm, with a constant independent of the cutoff.
Plancherel's identity gives
\[
0\le\frac{\nu}{(2\pi)^2}\int_a^b\int_{\R^2}
\ind_{\{|\xi|>N\}}|\wh\omega^\nu(t,\xi)|^2\dd\xi\dd t
\le\nu\int_a^b\|\omega^\nu(t)\|_{L^2(\R^2)}^2\dd t.
\]
Applying the preceding estimates to the right-hand side proves
both Fourier conclusions uniformly for $N>0$.
\end{proof}

\section{Proof of the main theorem}\label{sec:plateau-proof}

The balance \eqref{eq:high-balance} reduces the proof of the uniform
transfer limit to evaluating the endpoint energy increment and
showing that the viscous loss vanishes. We compare the increment
with the Gaussian energy. For $0<t<1$, we define
\begin{equation}\label{eq:Gaussian}
H_t^\nu(x)=\frac1{4\pi\nu t}e^{-|x|^2/(4\nu t)},
\qquad Q^\nu(t,x)=\sum_k\Gamma_kH_t^\nu(x-z_k(t)),
\end{equation}
and we set
\begin{equation}\label{eq:error}
r_\nu(t)=\|\omega^\nu(t)-Q^\nu(t)\|_{L^1(\R^2)}.
\end{equation}
We use the following approximation result.

\begin{proposition}[{\cite{NguyenCollapse}}]\label{prop:approximation}
For each fixed $0<\gamma<1$, there exists $\nu_\gamma>0$ such that
\begin{equation}\label{eq:approximation}
\sup_{0<t\le1-\nu^\gamma}r_\nu(t)
\lesssim_\gamma \nu^{1-\gamma},
\qquad 0<\nu\le\nu_\gamma.
\end{equation}
\end{proposition}

For every $s\in(1/2,1)$, vorticity perturbations of size
$O(\nu^s)$ in $L^1(\R^2)$, imposed at a fixed positive time
before collapse, can produce order-one $L^1(\R^2)$ separation near
collapse between Navier--Stokes solutions with the same viscosity
and total circulation~\cite{NguyenAmplification}.
Beyond collision, questions of continuation, nonuniqueness,
and selection by regularization arise~\cite{DM21,DMR24}.

We denote the velocity induced by the Gaussian approximation by
$u_H^\nu=K*Q^\nu$.

\begin{lemma}\label{lem:high-energy-error}
We fix $\nu,N>0$ and $0<t<1$.
We assume that $r_\nu(t)\le1$ and $\nu tN^2\le1$.
Then the energy difference satisfies
\begin{equation}\label{eq:high-energy-error}
\frac1{2(2\pi)^2}\left|\int_{\R^2}
\ind_{\{|\xi|>N\}}
\bigl(|\wh u^\nu(t)|^2-|\wh{u_H^\nu}(t)|^2\bigr)\dd\xi\right| \lesssim r_\nu(t)\left(1+\log\frac1{\nu tN^2r_\nu(t)}\right).
\end{equation}
We interpret the right-hand side as zero when $r_\nu(t)=0$.
\end{lemma}

\begin{proof}
The estimate \eqref{eq:high-energy-error} is immediate when
$r_\nu(t)=0$, because the two vorticities and their induced velocities
then agree. We therefore assume that $r_\nu(t)>0$, fix $t$, and omit
the time argument in Fourier transforms. To obtain the logarithmic
bound, it suffices to prove that, for every auxiliary cutoff $S\ge N$,
\[
\frac1{2(2\pi)^2}\left|\int_{\R^2} \ind_{\{|\xi|>N\}} \bigl(|\wh u^\nu|^2-|\wh{u_H^\nu}|^2\bigr)\dd\xi\right| \lesssim r_\nu(t)\log\frac SN+\frac1{\nu tS^2}.
\]
Indeed, the choice $S=(\nu tr_\nu(t))^{-1/2}$ makes the second term
of order $r_\nu(t)$, and the assumptions $r_\nu(t)\le1$ and
$\nu tN^2\le1$ ensure that this cutoff satisfies $S\ge N$.

To prove the auxiliary estimate, we separate the frequencies below
and above $S$. The Biot--Savart multiplier \eqref{eq:BS-multiplier}
gives the identities $\wh{u^\nu}(\xi) =-i\frac{\xi^\perp}{|\xi|^2}\wh\omega^\nu(\xi), \quad \wh{u_H^\nu}(\xi) =-i\frac{\xi^\perp}{|\xi|^2}\wh Q^\nu(\xi)$.
Since $|\xi^\perp|=|\xi|$, these identities imply that $|\wh{u^\nu}(\xi)|^2-|\wh{u_H^\nu}(\xi)|^2 =\frac{|\wh\omega^\nu(\xi)|^2-|\wh Q^\nu(\xi)|^2}{|\xi|^2}$.
The triangle inequality, followed by the splitting at $S$, therefore
yields
\[
\begin{aligned}
&\frac1{2(2\pi)^2}\left|\int_{\R^2}
\ind_{\{|\xi|>N\}}
\bigl(|\wh u^\nu|^2-|\wh{u_H^\nu}|^2\bigr)\dd\xi\right|\\
&\quad\lesssim
\int_{\R^2}\ind_{\{|\xi|>N\}}
\frac{\bigl||\wh\omega^\nu|^2-|\wh Q^\nu|^2\bigr|}
{|\xi|^2}\dd\xi\\
&\quad\le
\underbrace{\int_{\R^2}\ind_{\{N<|\xi|<S\}}
\frac{\bigl||\wh\omega^\nu|^2-|\wh Q^\nu|^2\bigr|}{|\xi|^2}
\dd\xi}_{(\mathrm{I})}
+\underbrace{\int_{\R^2}\ind_{\{|\xi|\ge S\}}
\frac{|\wh\omega^\nu|^2+|\wh Q^\nu|^2}{|\xi|^2}
\dd\xi}_{(\mathrm{II})}.
\end{aligned}
\]
Thus the auxiliary estimate follows from an $L^1$ comparison for
$(\mathrm{I})$ and $L^2$ bounds for $(\mathrm{II})$.

To bound $(\mathrm{I})$, it suffices to control the numerator by
$2\Gabs r_\nu(t)$. The factorization of the difference
of two squares and the reverse triangle inequality give
\[
\begin{aligned}
\bigl||\wh\omega^\nu(\xi)|^2-|\wh Q^\nu(\xi)|^2\bigr|
&=\bigl(|\wh\omega^\nu(\xi)|+|\wh Q^\nu(\xi)|\bigr)
\bigl||\wh\omega^\nu(\xi)|-|\wh Q^\nu(\xi)|\bigr|\\
&\le\bigl(|\wh\omega^\nu(\xi)|+|\wh Q^\nu(\xi)|\bigr)
|\wh\omega^\nu(\xi)-\wh Q^\nu(\xi)|.
\end{aligned}
\]
The definition of $r_\nu(t)$ bounds the last factor by
\[
\begin{aligned}
|\wh\omega^\nu(\xi)-\wh Q^\nu(\xi)|
&=\left|\int_{\R^2}e^{-ix\cdot\xi}
\bigl(\omega^\nu(t,x)-Q^\nu(t,x)\bigr)\dd x\right|\\
&\le\int_{\R^2}|\omega^\nu(t,x)-Q^\nu(t,x)|\dd x
=r_\nu(t).
\end{aligned}
\]
The circulation bound \eqref{eq:L1-bound} gives
$|\wh\omega^\nu(\xi)|\le\|\omega^\nu(t)\|_{L^1(\R^2)}\le\Gabs$.
The nonnegativity and unit mass of $H_t^\nu$ also give
$|\wh Q^\nu(\xi)|\le\|Q^\nu(t)\|_{L^1(\R^2)}
\le\sum_k|\Gamma_k|\int_{\R^2}H_t^\nu(x-z_k(t))\dd x
=\Gabs$.
Combining these bounds proves the required pointwise estimate, $\bigl||\wh\omega^\nu(\xi)|^2-|\wh Q^\nu(\xi)|^2\bigr| \le2\Gabs r_\nu(t)$.
The remaining radial integral is $\int_{\R^2}\ind_{\{N<|\xi|<S\}}\frac{\dd\xi}{|\xi|^2} =\int_0^{2\pi}\int_N^S\frac{\dd k}{k}\dd\theta =2\pi\log\frac SN$,
so the low-frequency contribution satisfies $(\mathrm{I})\lesssim \Gabs r_\nu(t)\log\frac SN$.

To estimate $(\mathrm{II})$, the inequality $|\xi|^{-2}\le S^{-2}$
on its domain and Plancherel's identity reduce the problem to
bounds for the two vorticities in $L^2$. The Gaussian definition gives
$\|H_t^\nu\|_{L^2(\R^2)}^2\lesssim(\nu t)^{-2}
\int_{\R^2}e^{-|x|^2/(2\nu t)}\dd x\lesssim(\nu t)^{-1}$.
The triangle inequality therefore yields
$\|Q^\nu(t)\|_{L^2(\R^2)}^2
\le(\sum_k|\Gamma_k|\|H_t^\nu\|_{L^2(\R^2)})^2
\lesssim\Gabs^2/(\nu t)$, while
\eqref{eq:dissipation-bound} gives
$\|\omega^\nu(t)\|_{L^2(\R^2)}^2
\lesssim\Gabs^2/(\nu t)$.
Consequently, the high-frequency contribution satisfies
\[
(\mathrm{II}) \lesssim S^{-2} \bigl(\|\omega^\nu(t)\|_{L^2(\R^2)}^2 +\|Q^\nu(t)\|_{L^2(\R^2)}^2\bigr) \lesssim\frac{\Gabs^2}{\nu tS^2}.
\]
The initial circulations are fixed, so the bounds for the two
contributions establish the auxiliary estimate for every $S\ge N$.
With the choice of $S$ made at the beginning of the proof, its
right-hand side satisfies
$r_\nu(t)\log\frac SN+\frac1{\nu tS^2}
=r_\nu(t)\left(1+\frac12\log\frac1{\nu tN^2r_\nu(t)}\right)
\le r_\nu(t)\left(1+\log\frac1{\nu tN^2r_\nu(t)}\right)$.
The logarithm is nonnegative under the hypotheses of the lemma.
This proves \eqref{eq:high-energy-error}.
\end{proof}

\begin{lemma}\label{lem:Gaussian-high}
The spectral profile $F$ is defined by \eqref{eq:profiles}.
For $\delta>0$ and $q\ge0$, we define
\begin{equation}\label{eq:interaction-integral}
I_\delta(q)=\frac1{4\pi}\int_q^\infty e^{-2\delta v^2}
\bigl(F(v)-\Gtot^2\bigr)\frac{\dd v}{v},
\qquad E_1(x)=\int_x^\infty e^{-s}\frac{\dd s}{s},\quad x>0.
\end{equation}
For every $\nu,N>0$ and $0<t<1$, the Gaussian energy satisfies
\begin{equation}\label{eq:Gaussian-high}
\frac1{2(2\pi)^2}\int_{\R^2}
\ind_{\{|\xi|>N\}}|\wh{u_H^\nu}(t)|^2\dd\xi =\frac{\Gtot^2}{8\pi}E_1(2\nu tN^2)
+I_{\nu t/(1-t)}\!\left(N\sqrt{1-t}\right).
\end{equation}
\end{lemma}

\begin{proof}
To prove \eqref{eq:Gaussian-high} for fixed $\nu,N>0$ and
$0<t<1$, it suffices to establish the radial representation
\begin{equation}\label{eq:Gaussian-high-radial}
\frac1{2(2\pi)^2}\int_{\R^2}
\ind_{\{|\xi|>N\}}|\wh{u_H^\nu}(t)|^2\dd\xi
=\frac1{4\pi}\int_N^\infty e^{-2\nu tk^2}F(k\sqrt{1-t})\frac{\dd k}{k}.
\end{equation}
Indeed, separating the constant part of the profile splits the
right-hand side into
\[
\begin{aligned}
&\frac1{4\pi}\int_N^\infty
e^{-2\nu tk^2}F(k\sqrt{1-t})\frac{\dd k}{k}\\
&\quad=\underbrace{\frac{\Gtot^2}{4\pi}
\int_N^\infty e^{-2\nu tk^2}\frac{\dd k}{k}}_{(\mathrm{I})}
+\underbrace{\frac1{4\pi}\int_N^\infty e^{-2\nu tk^2}
\bigl(F(k\sqrt{1-t})-\Gtot^2\bigr)
\frac{\dd k}{k}}_{(\mathrm{II})}.
\end{aligned}
\]
The substitution $s=2\nu tk^2$ evaluates the first term as $(\mathrm{I}) =\frac{\Gtot^2}{8\pi} \int_{2\nu tN^2}^\infty e^{-s}\frac{\dd s}{s} =\frac{\Gtot^2}{8\pi}E_1(2\nu tN^2)$.
To identify the second term, we set $q=N\sqrt{1-t},\quad \delta=\frac{\nu t}{1-t}$,
and change variables to $v=k\sqrt{1-t}$. The definition of
$I_\delta$ in \eqref{eq:interaction-integral} then gives $(\mathrm{II}) =\frac1{4\pi}\int_q^\infty e^{-2\delta v^2} \bigl(F(v)-\Gtot^2\bigr)\frac{\dd v}{v} =I_\delta(q)$.
Thus the two evaluations yield \eqref{eq:Gaussian-high} once
\eqref{eq:Gaussian-high-radial} has been proved.

It remains to establish that radial representation. The Fourier
transform of the Gaussian sum in \eqref{eq:Gaussian} is $\wh Q^\nu(t,\xi) =e^{-\nu t|\xi|^2}\sum_k\Gamma_k e^{-iz_k(t)\cdot\xi}$.
The trajectory formula \eqref{eq:centers} scales all positions by
$\sqrt{1-t}$ and applies a common rotation. This rotation only
shifts the angular integration variable, so
\eqref{eq:F-angular} gives $\frac1{2\pi}\int_0^{2\pi} \left|\sum_j\Gamma_j e^{-ik\mathbf e_\theta\cdot z_j(t)}\right|^2 \dd\theta=F(k\sqrt{1-t})$.
The Biot--Savart multiplier and polar coordinates therefore imply
\[
\begin{aligned}
&\frac1{2(2\pi)^2}\int_{\R^2}
\ind_{\{|\xi|>N\}}|\wh{u_H^\nu}(t)|^2\dd\xi\\
&\quad=\frac1{2(2\pi)^2}\int_N^\infty e^{-2\nu tk^2}
\left(\int_0^{2\pi}
\left|\sum_j\Gamma_j e^{-ik\mathbf e_\theta\cdot z_j(t)}\right|^2
\dd\theta\right)\frac{\dd k}{k}\\
&\quad=\frac1{4\pi}\int_N^\infty
e^{-2\nu tk^2}F(k\sqrt{1-t})\frac{\dd k}{k}.
\end{aligned}
\]
This proves \eqref{eq:Gaussian-high-radial} and completes the
reduction to \eqref{eq:Gaussian-high}.
\end{proof}

\begin{lemma}\label{lem:interaction-identity}
The functions $I_\delta$ and $E_1$ are defined by \eqref{eq:interaction-integral}.
For every $\delta>0$, we have
\begin{equation}\label{eq:full-cross}
I_\delta(0)=\cH
-\frac1{4\pi}\sum_{i<j}\Gamma_i\Gamma_j
E_1\!\left(\frac{r_{ij}(0)^2}{8\delta}\right).
\end{equation}
\end{lemma}

\begin{proof}
The identity \eqref{eq:full-cross} states that the sum of
$I_\delta(0)$ and the $E_1$ correction equals $\cH$.
It is therefore enough to show that this sum is constant in
$\delta>0$ and that its limit as $\delta\to\infty$ is
$\cH$. More precisely, our two goals are
\[
\begin{aligned}
\frac{\dd}{\dd\delta}\left[
I_\delta(0)+\frac1{4\pi}\sum_{i<j}\Gamma_i\Gamma_j
E_1\!\left(\frac{r_{ij}(0)^2}{8\delta}\right)\right]
&=0,\qquad \delta>0,\\
\lim_{\delta\to\infty}\left[
I_\delta(0)+\frac1{4\pi}\sum_{i<j}\Gamma_i\Gamma_j
E_1\!\left(\frac{r_{ij}(0)^2}{8\delta}\right)\right]
&=\cH.
\end{aligned}
\]

To prove the first identity, we use the definition of $E_1$ in
\eqref{eq:interaction-integral}, which gives $E_1'(x)=-e^{-x}/x$ for
$x>0$. The chain rule consequently yields
\[
\frac{\dd}{\dd\delta} E_1\!\left(\frac{r_{ij}(0)^2}{8\delta}\right) =-\frac{e^{-r_{ij}(0)^2/(8\delta)}}{r_{ij}(0)^2/(8\delta)} \left(-\frac{r_{ij}(0)^2}{8\delta^2}\right) =\frac1\delta e^{-r_{ij}(0)^2/(8\delta)}.
\]
Thus the vanishing derivative follows once we establish $\frac{\dd}{\dd\delta}I_\delta(0) =-\frac1{4\pi\delta}\sum_{i<j}\Gamma_i\Gamma_j e^{-r_{ij}(0)^2/(8\delta)}$.
To justify differentiation under the defining integral, we fix a
compact interval $0<\delta_0\le\delta\le\delta_1$.
The estimate \eqref{eq:F-origin} bounds the absolute value of the
differentiated integrand by $Cv^5$ near zero. The representation
\eqref{eq:profiles}, the identity \eqref{eq:circulation-square},
and the Hankel expansion \eqref{eq:Hankel} bound it by
$Cv^{1/2}e^{-2\delta_0v^2}$ at infinity. Both bounds are integrable
and uniform on this interval, so differentiation under the integral
is valid. Expanding the spectral profile then gives
\[
\begin{aligned}
\frac{\dd}{\dd\delta}I_\delta(0)
&=-\frac1{2\pi}\int_0^\infty v e^{-2\delta v^2}
\bigl(F(v)-\Gtot^2\bigr)\dd v\\
&=-\frac1\pi\sum_{i<j}\Gamma_i\Gamma_j
\int_0^\infty v e^{-2\delta v^2}J_0(r_{ij}(0)v)\dd v.
\end{aligned}
\]
The required derivative formula follows from the standard
Gaussian--Bessel identity~\cite{DLMF},
\[
\int_0^\infty v e^{-2\delta v^2}J_0(rv)\dd v
=\frac1{4\delta}e^{-r^2/(8\delta)},
\qquad r,\delta>0.
\]
Applying this identity with $r=r_{ij}(0)$ gives $\frac{\dd}{\dd\delta}I_\delta(0) =-\frac1{4\pi\delta}\sum_{i<j}\Gamma_i\Gamma_j e^{-r_{ij}(0)^2/(8\delta)}$.
This expression cancels the derivative of the $E_1$ correction, so the bracketed expression at the beginning
of the proof is constant on $(0,\infty)$.

To identify this constant, it suffices to prove the limits $I_\delta(0)\to 0, \quad \sum_{i<j}\Gamma_i\Gamma_j E_1\!\left(\frac{r_{ij}(0)^2}{8\delta}\right)\to 4\pi\cH \quad(\delta\to\infty)$.
For the first limit, we seek an integrable bound independent of
$\delta$ in the definition of $I_\delta(0)$. The quartic estimate
\eqref{eq:F-origin} provides the bound near zero, while
\eqref{eq:profiles}, \eqref{eq:circulation-square}, and
\eqref{eq:Hankel} provide the bound at infinity. Together they give
\[
\begin{aligned}
\int_0^\infty\frac{|F(v)-\Gtot^2|}{v}\dd v
&=\int_0^1\frac{|F(v)-\Gtot^2|}{v}\dd v
+\int_1^\infty\frac{|F(v)-\Gtot^2|}{v}\dd v\\
&\lesssim\int_0^1v^3\dd v
+\int_1^\infty v^{-3/2}\dd v<\infty.
\end{aligned}
\]
The Gaussian factor is bounded by one and tends to zero for every
$v>0$ as $\delta\to\infty$. Dominated convergence therefore proves
that $I_\delta(0)\to0$.

For the second limit, the arguments of $E_1$ tend to zero, so we
need its logarithmic expansion with a controlled remainder.
For $0<x\le1$, its defining integral gives
\[
E_1(x)+\log x =\int_x^1\frac{e^{-s}-1}{s}\dd s +\int_1^\infty\frac{e^{-s}}s\dd s =-\gamma_{\rm E}+\int_0^x\frac{1-e^{-s}}s\dd s,
\]
where $\gamma_{\rm E}$ is the Euler--Mascheroni constant.
The elementary bound $0\le1-e^{-s}\le s$ implies that $\bigl|E_1(x)+\log x+\gamma_{\rm E}\bigr| \le\int_0^x\dd s=x$.
We now set $c_{ij}=r_{ij}(0)^2/8$ and apply this estimate with
$x=c_{ij}/\delta$, which gives
\[
E_1(c_{ij}/\delta)
=\log(8\delta)-\gamma_{\rm E}-2\log r_{ij}(0)
+O(c_{ij}/\delta)
\qquad(\delta\to\infty).
\]
The circulation identity \eqref{eq:invariants} and the definition
\eqref{eq:transfer-constant} give the cancellations $\sum_{i<j}\Gamma_i\Gamma_j=0, \quad -2\sum_{i<j}\Gamma_i\Gamma_j\log r_{ij}(0)=4\pi\cH$.
After multiplication by $\Gamma_i\Gamma_j$ and summation over
$i<j$, the expansion consequently becomes
\[
\begin{aligned}
\sum_{i<j}\Gamma_i\Gamma_jE_1(c_{ij}/\delta)
&=(\log(8\delta)-\gamma_{\rm E})\sum_{i<j}\Gamma_i\Gamma_j
-2\sum_{i<j}\Gamma_i\Gamma_j\log r_{ij}(0)\\
&\quad+O\!\left(\frac1\delta
\sum_{i<j}|\Gamma_i\Gamma_j|c_{ij}\right)
=4\pi\cH+O(\delta^{-1}).
\end{aligned}
\]
The initial configuration is fixed, so the error tends to zero.
The two limits identify the constant bracketed expression as
$\cH$, and \eqref{eq:full-cross} follows.
\end{proof}

\begin{lemma}\label{lem:interaction-estimates}
For $0<q\le1$ and $0<\delta\le1$, the interaction integral
$I_\delta(q)$ defined in \eqref{eq:interaction-integral} has the small-$q$ expansion
\begin{equation}\label{eq:cross-low}
I_\delta(q)=\cH+O(q^4+\delta e^{-c/\delta}),
\qquad c=\frac{r_{\min}(0)^2}{8}.
\end{equation}
For $q\ge1$ and every $\delta>0$, it satisfies the large-$q$ bound
\begin{equation}\label{eq:cross-tail}
|I_\delta(q)|\lesssim q^{-3/2}
\sum_{i<j}|\Gamma_i\Gamma_j|r_{ij}(0)^{-3/2}.
\end{equation}
The bound \eqref{eq:cross-tail} has a universal implicit constant
and holds for every $q>0$.
\end{lemma}

\begin{proof}
To prove the small-argument estimate \eqref{eq:cross-low}, it
suffices to control the change from $q=0$ and the correction to
$\cH$ at $q=0$ separately. Indeed, the triangle inequality
gives $|I_\delta(q)-\cH| \le |I_\delta(q)-I_\delta(0)| +|I_\delta(0)-\cH|$.
We therefore seek a bound of order $q^4$ for the first difference
and a bound of order $\delta e^{-c/\delta}$ for the second.

For the first difference, the definition
\eqref{eq:interaction-integral} reduces the estimate to the
spectral profile near zero, since $I_\delta(q)-I_\delta(0) =-\frac1{4\pi}\int_0^q e^{-2\delta v^2} \bigl(F(v)-\Gtot^2\bigr)\frac{\dd v}{v}$.
The quartic bound \eqref{eq:F-origin} gives $|F(v)-\Gtot^2| \lesssim v^4,\quad 0<v\le1$.
Since the Gaussian factor is at most one, this bound yields
\[
|I_\delta(q)-I_\delta(0)| \le\frac1{4\pi}\int_0^q \frac{|F(v)-\Gtot^2|}{v}\dd v \lesssim\int_0^q v^3\dd v =\frac{q^4}{4},\qquad 0<q\le1.
\]

For the second difference, the identity \eqref{eq:full-cross}
reduces the estimate to the exponential integral, because $I_\delta(0)-\cH =-\frac1{4\pi}\sum_{i<j}\Gamma_i\Gamma_j E_1\!\left(\frac{r_{ij}(0)^2}{8\delta}\right)$.
The definition of $E_1$ implies the bound $0<E_1(x)=\int_x^\infty\frac{e^{-s}}s\dd s \le\frac1x\int_x^\infty e^{-s}\dd s =\frac{e^{-x}}x, \quad x>0$.
We set $c_{ij}=r_{ij}(0)^2/8$, so that
$c=\min_{i<j}c_{ij}=r_{\min}(0)^2/8>0$.
Applying the bound with $x=c_{ij}/\delta$ and using
$c_{ij}\ge c$, we obtain
\[
\begin{aligned}
|I_\delta(0)-\cH|
&\le\frac1{4\pi}\sum_{i<j}|\Gamma_i\Gamma_j|E_1(c_{ij}/\delta)
\le\frac{\delta}{4\pi}\sum_{i<j}
\frac{|\Gamma_i\Gamma_j|}{c_{ij}}e^{-c_{ij}/\delta}\\
&\le\left(\frac1{4\pi}\sum_{i<j}
\frac{|\Gamma_i\Gamma_j|}{c_{ij}}\right)
\delta e^{-c/\delta}.
\end{aligned}
\]
The coefficient depends only on the fixed initial configuration.
The two estimates for the differences give $|I_\delta(q)-\cH| \lesssim q^4+\delta e^{-c/\delta}, \quad 0<q,\delta\le1$,
which proves \eqref{eq:cross-low}.

To prove the tail estimate \eqref{eq:cross-tail} uniformly in
$\delta$, we remove the Gaussian weight by integration by parts
and estimate the resulting undamped tail. For this purpose, we
define $B(v)=\frac1{4\pi}\int_v^\infty \bigl(F(s)-\Gtot^2\bigr)\frac{\dd s}{s}, \quad v>0$.
The representation \eqref{eq:profiles}, the identity
\eqref{eq:circulation-square}, and the Hankel expansion
\eqref{eq:Hankel} show that the integrand is $O(s^{-3/2})$ at
infinity. This integral therefore converges absolutely, and
$B(v)\to0$ as $v\to\infty$.
We shall reduce the desired tail estimate to the bound $|B(v)|\le Cv^{-3/2} \sum_{i<j}|\Gamma_i\Gamma_j|r_{ij}(0)^{-3/2}, \quad v>0$,
with a universal constant $C$.

To justify this reduction, we differentiate the undamped tail
and obtain $B'(v)=-\frac{F(v)-\Gtot^2}{4\pi v}$.
Substitution into the definition gives
$I_\delta(q)=-\int_q^\infty e^{-2\delta v^2}B'(v)\dd v$,
and integration by parts yields
\[
\begin{aligned}
I_\delta(q)
&=-\left[e^{-2\delta v^2}B(v)\right]_q^\infty
-4\delta\int_q^\infty v e^{-2\delta v^2}B(v)\dd v\\
&=e^{-2\delta q^2}B(q)
-4\delta\int_q^\infty v e^{-2\delta v^2}B(v)\dd v.
\end{aligned}
\]
The boundary term at infinity vanishes because $B(v)\to0$.
The total mass of the differentiated Gaussian weight is $4\delta\int_q^\infty v e^{-2\delta v^2}\dd v =\left[-e^{-2\delta v^2}\right]_q^\infty =e^{-2\delta q^2}$.
Taking absolute values in the integration-by-parts formula
therefore yields $|I_\delta(q)| \le2e^{-2\delta q^2}\sup_{v\ge q}|B(v)|$.
The claimed estimate for $B(v)$ would now imply
\eqref{eq:cross-tail}, because $v^{-3/2}\le q^{-3/2}$ for
$v\ge q$ and $e^{-2\delta q^2}\le1$.

To bound $B(v)$, we use the identity
$F(s)-\Gtot^2=2\sum_{i<j}\Gamma_i\Gamma_jJ_0(r_{ij}(0)s)$,
which follows from \eqref{eq:profiles} and \eqref{eq:circulation-square}.
Changing variables in each integral gives $B(v)=\frac1{2\pi}\sum_{i<j}\Gamma_i\Gamma_j \int_{r_{ij}(0)v}^\infty\frac{J_0(s)}s\dd s$.
Thus it suffices to prove the scalar bound $\left|\int_y^\infty\frac{J_0(s)}s\dd s\right| \le Cy^{-3/2},\quad y>0$,
with a universal constant $C$.
For $y\ge1$, the case $m=1$ of \eqref{eq:app-Bessel-tails}
provides this estimate directly. For $0<y\le1$, the angular
representation \eqref{eq:app-bessel-angular} gives
$|J_0(s)|\le1$, so splitting the integral at $1$ yields
\[
\begin{aligned}
\left|\int_y^\infty\frac{J_0(s)}s\dd s\right|
&\le\int_y^1\frac{|J_0(s)|}{s}\dd s
+\left|\int_1^\infty\frac{J_0(s)}s\dd s\right|\\
&\le |\log y|+C
\lesssim y^{-3/2}.
\end{aligned}
\]
The last comparison has a universal constant on $(0,1]$.
The scalar estimate consequently holds for every $y>0$.
Applying it with $y=r_{ij}(0)v$ in the pair representation gives
\[
|B(v)| \le\frac1{2\pi}\sum_{i<j}|\Gamma_i\Gamma_j| \left|\int_{r_{ij}(0)v}^\infty\frac{J_0(s)}s\dd s\right| \lesssim v^{-3/2} \sum_{i<j}|\Gamma_i\Gamma_j|r_{ij}(0)^{-3/2}.
\]
Returning to the reduction for $I_\delta(q)$, we obtain
\[
\begin{aligned}
|I_\delta(q)|
&\le2e^{-2\delta q^2}\sup_{v\ge q}|B(v)|
\lesssim e^{-2\delta q^2}q^{-3/2}
\sum_{i<j}|\Gamma_i\Gamma_j|r_{ij}(0)^{-3/2}\\
&\lesssim q^{-3/2}
\sum_{i<j}|\Gamma_i\Gamma_j|r_{ij}(0)^{-3/2}.
\end{aligned}
\]
All constants in this tail estimate are universal, and the
argument applies to every $q,\delta>0$.
This proves \eqref{eq:cross-tail}.
\end{proof}

\begin{proof}[Proof of the main theorem]
We prove the two estimates in \eqref{eq:plateau-estimates} by reducing the nonlinear-transfer error to the energy error and the viscous loss. The integrated balance \eqref{eq:high-balance} gives the exact relation
\[
\begin{aligned}
&-\frac1{(2\pi)^2}\int_{A_\nu}^{B_\nu}\Real
\int_{\R^2}\ind_{\{|\xi|>N\}}
\ol{\wh u^\nu(t)}\cdot
\wh{(u^\nu\cdot\nabla)u^\nu}(t)\dd\xi\dd t-\cH\\
&\quad=\left[\frac1{2(2\pi)^2}\int_{\R^2}
\ind_{\{|\xi|>N\}}
\Bigl(|\wh u^\nu(B_\nu)|^2-|\wh u^\nu(A_\nu)|^2\Bigr)
\dd\xi-\cH\right]\\
&\qquad+\frac{\nu}{(2\pi)^2}\int_{A_\nu}^{B_\nu}
\int_{\R^2}\ind_{\{|\xi|>N\}}
|\wh\omega^\nu(t)|^2\dd\xi\dd t.
\end{aligned}
\]
The desired estimates therefore follow if the bracketed energy error satisfies the claimed bound and the viscous loss is at most $C\nu^{\gamma_-}$, with both bounds uniform for $\nu^{-\alpha_-/2}\le N\le\nu^{-\alpha_+/2}$. The energy error will be reduced to the Gaussian increment and two endpoint comparison errors after we establish the required bound for the viscous loss.

The bound for the viscous loss reduces to a bound for $\log(B_\nu/A_\nu)$, because Plancherel's identity and \eqref{eq:dissipation-bound} imply
\[
\begin{aligned}
0&\le\frac{\nu}{(2\pi)^2}\int_{A_\nu}^{B_\nu}
\int_{\R^2}\ind_{\{|\xi|>N\}}
|\wh\omega^\nu(t)|^2\dd\xi\dd t\\
&\le\nu\int_{A_\nu}^{B_\nu}
\|\omega^\nu(t)\|_{L^2(\R^2)}^2\dd t
\lesssim\Gabs^2\int_{A_\nu}^{B_\nu}\frac{\dd t}{t}
=\Gabs^2\log\frac{B_\nu}{A_\nu}.
\end{aligned}
\]
The required logarithmic bound follows from the length of the observation interval. For sufficiently small $\nu$, the observation times satisfy \(0<B_\nu-A_\nu=\nu^{\gamma_-}-\nu^{\gamma_+}\le\nu^{\gamma_-}\) and \(A_\nu\ge1/2\). These inequalities and the bound $\log(1+x)\le x$ for $x\ge0$ give \(0<\log(B_\nu/A_\nu)=\log(1+(B_\nu-A_\nu)/A_\nu)\le(B_\nu-A_\nu)/A_\nu\le2\nu^{\gamma_-}\). This proves the required $O(\nu^{\gamma_-})$ bound for the viscous loss, uniformly in $N$. It remains to estimate the energy error.

The energy error can be estimated once we control the Gaussian energy increment and the comparison error at each endpoint. To this end, we use $u_H^\nu=K*Q^\nu$, where the Gaussian sum in \eqref{eq:Gaussian} is \(Q^\nu(t,x)=\sum_{k=1}^3\frac{\Gamma_k}{4\pi\nu t}\exp(-|x-z_k(t)|^2/(4\nu t))\). We add and subtract the Gaussian energy at each observation time to obtain
\[
\frac1{2(2\pi)^2}\int_{\R^2}
\ind_{\{|\xi|>N\}}
\Bigl(|\wh u^\nu(B_\nu)|^2
-|\wh u^\nu(A_\nu)|^2\Bigr)\dd\xi
=(\mathrm{I})+(\mathrm{II})-(\mathrm{III}),
\]
where the three terms are defined by
\[
\begin{aligned}
(\mathrm{I})
&=\frac1{2(2\pi)^2}\int_{\R^2}
\ind_{\{|\xi|>N\}}
\Bigl(|\wh u_H^\nu(B_\nu)|^2
-|\wh u_H^\nu(A_\nu)|^2\Bigr)\dd\xi,\\
(\mathrm{II})
&=\frac1{2(2\pi)^2}\int_{\R^2}
\ind_{\{|\xi|>N\}}
\Bigl(|\wh u^\nu(B_\nu)|^2
-|\wh u_H^\nu(B_\nu)|^2\Bigr)\dd\xi,\\
(\mathrm{III})
&=\frac1{2(2\pi)^2}\int_{\R^2}
\ind_{\{|\xi|>N\}}
\Bigl(|\wh u^\nu(A_\nu)|^2
-|\wh u_H^\nu(A_\nu)|^2\Bigr)\dd\xi.
\end{aligned}
\]
The triangle inequality now reduces the energy error to three terms:
\[
\begin{aligned}
&\left|\frac1{2(2\pi)^2}\int_{\R^2}
\ind_{\{|\xi|>N\}}
\Bigl(|\wh u^\nu(B_\nu)|^2
-|\wh u^\nu(A_\nu)|^2\Bigr)\dd\xi-\cH\right|\\
&\quad\le |(\mathrm{I})-\cH|
+|(\mathrm{II})|+|(\mathrm{III})|.
\end{aligned}
\]
The required bound will follow from an estimate for the Gaussian error $|(\mathrm{I})-\cH|$ and logarithmic comparison estimates for $|(\mathrm{II})|$ and $|(\mathrm{III})|$.

The Gaussian error reduces further to an $E_1$-difference and two interaction terms through \eqref{eq:Gaussian-high}. The functions in that formula are defined in \eqref{eq:interaction-integral} by \(I_\delta(q)=\frac1{4\pi}\int_q^\infty e^{-2\delta v^2}(F(v)-\Gtot^2)\,\dd v/v\) and \(E_1(x)=\int_x^\infty e^{-s}\,\dd s/s\). The difference of the two instances of \eqref{eq:Gaussian-high}, at $B_\nu$ and $A_\nu$, gives
\[
\begin{aligned}
(\mathrm{I})
&=\frac{\Gtot^2}{8\pi}
\bigl[E_1(2\nu B_\nu N^2)-E_1(2\nu A_\nu N^2)\bigr]\\
&\qquad+I_{\nu B_\nu/(1-B_\nu)}\!\left(N\sqrt{1-B_\nu}\right)
-I_{\nu A_\nu/(1-A_\nu)}\!\left(N\sqrt{1-A_\nu}\right).
\end{aligned}
\]
We subtract $\cH$ and apply the triangle inequality to obtain
\[
\begin{aligned}
|(\mathrm{I})-\cH|
&\le\frac{\Gtot^2}{8\pi}
\left|E_1(2\nu B_\nu N^2)-E_1(2\nu A_\nu N^2)\right|\\
&\quad+\left|I_{\nu A_\nu/(1-A_\nu)}
\!\left(N\sqrt{1-A_\nu}\right)\right|\\
&\quad+\left|I_{\nu B_\nu/(1-B_\nu)}
\!\left(N\sqrt{1-B_\nu}\right)-\cH\right|.
\end{aligned}
\]
The first contribution will be bounded by the observation interval length. The other two contributions will be controlled by the large- and small-argument estimates for $I_\delta$, respectively.

The required $O(\nu^{\gamma_-})$ estimate for the $E_1$-difference follows if the absolute value of its time derivative is bounded by $1/t$. The definition of $E_1$ gives the derivative \(\frac{\dd}{\dd t}E_1(2\nu tN^2)=-\frac{e^{-2\nu tN^2}}{2\nu tN^2}\,2\nu N^2=-e^{-2\nu tN^2}/t\). We integrate the derivative over $[A_\nu,B_\nu]$ and use the logarithmic bound proved above to obtain
\[
\begin{aligned}
&\frac{\Gtot^2}{8\pi}
\left|E_1(2\nu B_\nu N^2)-E_1(2\nu A_\nu N^2)\right|\\
&\quad=\frac{\Gtot^2}{8\pi}
\int_{A_\nu}^{B_\nu}\frac{e^{-2\nu tN^2}}t\dd t
\le\frac{\Gtot^2}{8\pi}
\log\frac{B_\nu}{A_\nu}
\lesssim\nu^{\gamma_-}\longrightarrow0.
\end{aligned}
\]

The required estimate for the interaction term at $A_\nu$ follows from \eqref{eq:cross-tail} once its argument is bounded from below in terms of $\nu$. The lower cutoff bound yields \(N\sqrt{1-A_\nu}=N\nu^{\gamma_-/2}\ge\nu^{-(\alpha_--\gamma_-)/2}\longrightarrow\infty\). We apply the uniform tail estimate \eqref{eq:cross-tail} to this argument and obtain
\[
\begin{aligned}
\left|I_{\nu A_\nu/(1-A_\nu)}\!\left(N\sqrt{1-A_\nu}\right)\right|
&\lesssim\bigl(N\sqrt{1-A_\nu}\bigr)^{-3/2}
\sum_{i<j}|\Gamma_i\Gamma_j|r_{ij}(0)^{-3/2}\\
&\lesssim\nu^{3(\alpha_--\gamma_-)/4}\longrightarrow0.
\end{aligned}
\]
The initial configuration is fixed independently of $\nu$, so the distance sum is absorbed into the implicit constant.

The remaining interaction error, at $B_\nu$, is controlled by \eqref{eq:cross-low} provided that both arguments of $I_\delta(q)$ lie in $(0,1]$. The upper cutoff bound and the definition of $B_\nu$ verify this condition for small $\nu$, since \(N\sqrt{1-B_\nu}=N\nu^{\gamma_+/2}\le\nu^{(\gamma_+-\alpha_+)/2}\longrightarrow0\) and \(0<\nu B_\nu/(1-B_\nu)=B_\nu\nu^{1-\gamma_+}\le\nu^{1-\gamma_+}\longrightarrow0\). These bounds hold uniformly throughout the prescribed cutoff range. We substitute the two arguments into \eqref{eq:cross-low} and obtain
\[
\begin{aligned}
&\left|I_{\nu B_\nu/(1-B_\nu)}
\!\left(N\sqrt{1-B_\nu}\right)-\cH\right|\\
&\quad\lesssim\bigl(N\sqrt{1-B_\nu}\bigr)^4
+\frac{\nu B_\nu}{1-B_\nu}
\exp\!\left(-\frac{c(1-B_\nu)}{\nu B_\nu}\right)\\
&\quad\le\nu^{2(\gamma_+-\alpha_+)}
+\nu^{1-\gamma_+}e^{-c/\nu^{1-\gamma_+}}
\longrightarrow0.
\end{aligned}
\]
The bounds for the $E_1$-difference and the two interaction terms now give the required estimate for the Gaussian error:
\[
|(\mathrm{I})-\cH|\lesssim\nu^{\gamma_-}+\nu^{2(\gamma_+-\alpha_+)}
+\nu^{3(\alpha_--\gamma_-)/4}+\nu^{1-\gamma_+}e^{-c/\nu^{1-\gamma_+}}.
\]

The two remaining comparison terms will satisfy the required $O(\nu^{1-\gamma_+}(1+|\log\nu|))$ bound if we can apply \eqref{eq:high-energy-error} and control its logarithmic factor. We verify the hypotheses at both observation times. Proposition~\ref{prop:approximation}, applied with $\gamma=\gamma_+$, provides a constant $C_0\ge1$ such that, for sufficiently small $\nu$, \(0\le r_\nu(t)=\|\omega^\nu(t)-Q^\nu(t)\|_{L^1(\R^2)}\le C_0\nu^{1-\gamma_+}\le1\) for \(t\in\{A_\nu,B_\nu\}\). The observation times satisfy $t\ge1/2$, so the prescribed cutoff range also gives \(\tfrac12\nu^{1-\alpha_-}\le\nu tN^2\le\nu^{1-\alpha_+}\le1\). The logarithmic comparison bound can be evaluated at the upper bound for $r_\nu(t)$ because its right-hand side is increasing in $r$. Indeed, the preceding bound $\nu tN^2\le1$ implies \(\frac{\dd}{\dd r}[r(1+\log(1/(\nu tN^2r)))]=\log(1/(\nu tN^2r))\ge0\) for \(0<r\le1\). The function $r\mapsto r(1+\log(1/(\nu tN^2r)))$ extends continuously with value zero at $r=0$. We use that value when $r_\nu(t)=0$, in which case the two vorticities agree and the corresponding comparison term vanishes.

We first apply \eqref{eq:high-energy-error} at $t=B_\nu$. The monotonicity just established and the bounds on $r_\nu(B_\nu)$ and $\nu B_\nu N^2$ yield
\[
\begin{aligned}
|(\mathrm{II})|&\lesssim r_\nu(B_\nu)
\left(1+\log\frac1{\nu B_\nu N^2r_\nu(B_\nu)}\right)\\
&\le C_0\nu^{1-\gamma_+}
\left(1+\log\frac{2}{C_0\nu^{2-\alpha_--\gamma_+}}\right)
\lesssim\nu^{1-\gamma_+}(1+|\log\nu|).
\end{aligned}
\]
We then apply \eqref{eq:high-energy-error} at $t=A_\nu$. The same bounds give
\[
\begin{aligned}
|(\mathrm{III})|&\lesssim r_\nu(A_\nu)
\left(1+\log\frac1{\nu A_\nu N^2r_\nu(A_\nu)}\right)\\
&\le C_0\nu^{1-\gamma_+}
\left(1+\log\frac{2}{C_0\nu^{2-\alpha_--\gamma_+}}\right)
\lesssim\nu^{1-\gamma_+}(1+|\log\nu|).
\end{aligned}
\]

Combining these three estimates and bounding the exponentially small term by $\nu^{1-\gamma_+}$, we obtain
\begin{equation}\label{eq:endpoint-error}
\begin{aligned}
&\sup_{\substack{\nu^{-\alpha_-/2}\le N\\N\le\nu^{-\alpha_+/2}}}
\left|\frac1{2(2\pi)^2}\int_{\R^2}
\ind_{\{|\xi|>N\}}
\Bigl(|\wh u^\nu(B_\nu)|^2-|\wh u^\nu(A_\nu)|^2\Bigr)
\dd\xi-\cH\right|\\
&\quad\le\sup_{\substack{\nu^{-\alpha_-/2}\le N\\N\le\nu^{-\alpha_+/2}}}
\bigl(|(\mathrm{I})-\cH|+|(\mathrm{II})|+|(\mathrm{III})|\bigr)\\
&\quad\lesssim\nu^{1-\gamma_+}(1+|\log\nu|)
+\nu^{\gamma_-}+\nu^{2(\gamma_+-\alpha_+)}
+\nu^{3(\alpha_--\gamma_-)/4}.
\end{aligned}
\end{equation}
This is the energy estimate in \eqref{eq:plateau-estimates}. The integrated energy balance, together with \eqref{eq:endpoint-error} and the $O(\nu^{\gamma_-})$ viscous-loss bound, now gives the nonlinear-transfer estimate with the same rate after an increase of $C$. The ordering $0<\gamma_-<\alpha_-<\alpha_+<\gamma_+<1$ makes every exponent positive, so both errors tend to zero uniformly over the stated cutoff range.

The dyadic assertions follow from the uniform estimates once we verify that all the cutoffs in \eqref{eq:dyadic-cutoffs} lie in the prescribed range and that their number tends to infinity. For every integer $0\le j\le J_\nu$, the floor in the definition of $J_\nu$ gives
\[
\nu^{-\alpha_-/2}\le N_{j,\nu}=2^j\nu^{-\alpha_-/2}
\le2^{J_\nu}\nu^{-\alpha_-/2}
\le\nu^{-(\alpha_+-\alpha_-)/2}\nu^{-\alpha_-/2}
=\nu^{-\alpha_+/2}.
\]
Thus every dyadic cutoff is admissible, and the two error bounds hold simultaneously for all $j$. The same definition and the strict inequality $\alpha_+>\alpha_-$ give \(J_\nu\ge\frac{\alpha_+-\alpha_-}{2\log2}|\log\nu|-1\longrightarrow\infty\). The strict negativity of $\cH$ therefore identifies a net transfer of $-\cH+o(1)$ toward lower frequencies across each of the $J_\nu+1$ cutoffs. The error is uniform in $j$, which completes the dyadic assertions and the proof.
\end{proof}

\begin{proof}[Proof of Corollary~\ref{cor:energy-gain}]
We reduce the two band limits in \eqref{eq:band-gain} to vanishing energy and transfer increments above the fixed cutoff $k_0$. The main theorem already gives the limiting increments above $N$, and the energy increment and integrated transfer in the band are the corresponding differences because \(\ind_{\{k_0<|\xi|\le N\}}=\ind_{\{|\xi|>k_0\}}-\ind_{\{|\xi|>N\}}\) for \(N>k_0\). It therefore suffices to show that both quantities above $k_0$ tend to zero. The transfer bound will follow from the energy bound and the vanishing viscous loss through \eqref{eq:high-balance}.

The fixed-cutoff energy bound reduces to the Gaussian increment and two endpoint comparison errors. We add and subtract the Gaussian energy at each endpoint and use the triangle inequality to obtain
\[
\begin{aligned}
&\left|\frac1{2(2\pi)^2}\int_{\R^2}
\ind_{\{|\xi|>k_0\}}
\bigl[|\wh u^\nu(t)|^2\bigr]_{t=A_\nu}^{t=B_\nu}\dd\xi\right|\\
&\quad\le
\left|\frac1{2(2\pi)^2}\int_{\R^2}
\ind_{\{|\xi|>k_0\}}
\bigl[|\wh{u_H^\nu}(t)|^2\bigr]_{t=A_\nu}^{t=B_\nu}\dd\xi\right|\\
&\qquad+
\sum_{t\in\{A_\nu,B_\nu\}}\frac1{2(2\pi)^2}
\left|\int_{\R^2}\ind_{\{|\xi|>k_0\}}
\bigl(|\wh u^\nu(t)|^2
-|\wh{u_H^\nu}(t)|^2\bigr)\dd\xi\right|.
\end{aligned}
\]
The two comparison errors are controlled by \eqref{eq:high-energy-error} once its hypotheses hold with $N=k_0$. At both endpoints, the observation times and \eqref{eq:approximation} give, for sufficiently small $\nu$, \(1/2\le t<1\), \(r_\nu(t)\le C_0\nu^{1-\gamma_+}\le1\), and \(\tfrac12\nu k_0^2\le\nu tk_0^2\le1\). We apply \eqref{eq:high-energy-error} and then evaluate its right-hand side at the upper bound for $r_\nu(t)$. This step uses the monotonicity of $r\mapsto r(1+\log(1/(\nu tk_0^2r)))$ on $0<r\le1$, which was proved in the argument for the main theorem. The resulting estimate is
\[
\begin{aligned}
&\frac1{2(2\pi)^2}\left|\int_{\R^2}
\ind_{\{|\xi|>k_0\}}
\bigl(|\wh u^\nu(t)|^2-|\wh{u_H^\nu}(t)|^2\bigr)\dd\xi\right|\\
&\quad\lesssim r_\nu(t)
\left(1+\log\frac1{\nu tk_0^2r_\nu(t)}\right)\\
&\quad\le C_0\nu^{1-\gamma_+}
\left(1+\log\frac2{C_0k_0^2\nu^{2-\gamma_+}}\right)
\lesssim_{k_0}\nu^{1-\gamma_+}(1+|\log\nu|),
\qquad t\in\{A_\nu,B_\nu\}.
\end{aligned}
\]
The estimate also holds at $r_\nu(t)=0$ by the continuous interpretation of its right-hand side, because the comparison error then vanishes.

It remains to estimate the Gaussian energy increment. Formula \eqref{eq:Gaussian-high} expresses it as an $E_1$-difference and two deviations from $\cH$:
\[
\begin{aligned}
&\frac1{2(2\pi)^2}\int_{\R^2}\ind_{\{|\xi|>k_0\}}
\Bigl(|\wh{u_H^\nu}(B_\nu)|^2
-|\wh{u_H^\nu}(A_\nu)|^2\Bigr)\dd\xi\\
&\quad=\frac{\Gtot^2}{8\pi}
\Bigl[E_1(2\nu B_\nu k_0^2)-E_1(2\nu A_\nu k_0^2)\Bigr]\\
&\qquad+
\left[I_{\nu B_\nu/(1-B_\nu)}
\!\left(k_0\sqrt{1-B_\nu}\right)-\cH\right]
-\left[I_{\nu A_\nu/(1-A_\nu)}
\!\left(k_0\sqrt{1-A_\nu}\right)-\cH\right].
\end{aligned}
\]
The $E_1$-difference has the required $O(\nu^{\gamma_-})$ bound by the derivative calculation in the proof of the main theorem. More precisely, that calculation gives
\[
\left|E_1(2\nu B_\nu k_0^2)-E_1(2\nu A_\nu k_0^2)\right|
=\int_{A_\nu}^{B_\nu}\frac{e^{-2\nu tk_0^2}}t\dd t
\le\log\frac{B_\nu}{A_\nu}\lesssim\nu^{\gamma_-}.
\]
The two interaction deviations are controlled by \eqref{eq:cross-low} once both of its arguments are small. At both endpoints, the definitions of the observation times give \(k_0\sqrt{1-t}\le k_0\nu^{\gamma_-/2}\longrightarrow0\) and \(\nu t/(1-t)\le\nu^{1-\gamma_+}\longrightarrow0\). These inequalities verify the required hypotheses for sufficiently small $\nu$. We therefore apply \eqref{eq:cross-low} at each endpoint and obtain
\[
\left|I_{\nu t/(1-t)}\!\left(k_0\sqrt{1-t}\right)-\cH\right|
\lesssim_{k_0}\nu^{2\gamma_-}+\nu^{1-\gamma_+}e^{-c/\nu^{1-\gamma_+}},
\qquad t\in\{A_\nu,B_\nu\}.
\]
Combining these three estimates, we absorb $\nu^{2\gamma_-}$ into $\nu^{\gamma_-}$ and the exponential term into $\nu^{1-\gamma_+}$ to obtain
\begin{equation}\label{eq:fixed-cutoff-energy}
\begin{aligned}
&\left|\frac1{2(2\pi)^2}\int_{\R^2}
\ind_{\{|\xi|>k_0\}}
\bigl[|\wh u^\nu(t)|^2\bigr]_{t=A_\nu}^{t=B_\nu}\dd\xi\right|\\
&\qquad\lesssim_{k_0}
\nu^{1-\gamma_+}(1+|\log\nu|)+\nu^{\gamma_-}
\longrightarrow0.
\end{aligned}
\end{equation}

The fixed-cutoff transfer will now vanish if the viscous loss vanishes. Indeed, the time integral of \eqref{eq:high-balance} with $N=k_0$ gives
\[
\begin{aligned}
&-\frac1{(2\pi)^2}\int_{A_\nu}^{B_\nu}\Real
\int_{\R^2}\ind_{\{|\xi|>k_0\}}
\ol{\wh u^\nu(t)}\cdot
\wh{(u^\nu\cdot\nabla)u^\nu}(t)\dd\xi\dd t\\
&\quad=\frac1{2(2\pi)^2}\int_{\R^2}
\ind_{\{|\xi|>k_0\}}
\bigl[|\wh u^\nu(t)|^2\bigr]_{t=A_\nu}^{t=B_\nu}\dd\xi
+\frac{\nu}{(2\pi)^2}\int_{A_\nu}^{B_\nu}
\int_{\R^2}\ind_{\{|\xi|>k_0\}}
|\wh\omega^\nu(t)|^2\dd\xi\dd t.
\end{aligned}
\]
We control the remaining viscous term by \eqref{eq:integrated-dissipation-bound} and the observation interval estimate from the main proof:
\[
\frac{\nu}{(2\pi)^2}\int_{A_\nu}^{B_\nu}
\int_{\R^2}\ind_{\{|\xi|>k_0\}}
|\wh\omega^\nu(t)|^2\dd\xi\dd t
\le\nu\int_{A_\nu}^{B_\nu}\|\omega^\nu(t)\|_2^2\dd t
\overset{\eqref{eq:integrated-dissipation-bound}}{\lesssim}
\log\frac{B_\nu}{A_\nu}\lesssim\nu^{\gamma_-}.
\]
The integrated balance and \eqref{eq:fixed-cutoff-energy} therefore give the required transfer bound:
\begin{equation}\label{eq:fixed-cutoff-transfer}
\begin{aligned}
&\left|-\frac1{(2\pi)^2}\int_{A_\nu}^{B_\nu}\Real
\int_{\R^2}\ind_{\{|\xi|>k_0\}}
\ol{\wh u^\nu(t)}\cdot
\wh{(u^\nu\cdot\nabla)u^\nu}(t)\dd\xi\dd t\right|\\
&\qquad\lesssim_{k_0}
\nu^{1-\gamma_+}(1+|\log\nu|)+\nu^{\gamma_-}
\longrightarrow0.
\end{aligned}
\end{equation}

Every admissible cutoff exceeds $k_0$ for sufficiently small $\nu$. By the indicator identity at the start of the proof, each quantity in \eqref{eq:band-gain} is its counterpart above $k_0$ minus its counterpart above $N$. The bounds \eqref{eq:fixed-cutoff-energy} and \eqref{eq:fixed-cutoff-transfer}, together with the main theorem, therefore give the common limit $0-\cH=-\cH$, uniformly in $N$.
\end{proof}

\section*{Acknowledgments}

The author is supported by the NSFC Excellent Young Scientists Fund
(Overseas) and the Research Start-up Fund of Shanghai Jiao Tong University.

\subsection*{Data availability}

No datasets were generated or analysed during the current study.

\subsection*{Competing interests}

The author declares no competing interests.

\end{document}